\documentclass[12pt]{article}

\usepackage[authoryear]{natbib}
\usepackage{braket}
\usepackage{enumitem}

\usepackage{etoolbox}
\RequirePackage[OT1]{fontenc}
\usepackage{flafter}
\usepackage{makecell}
\usepackage{graphicx}
\usepackage[space]{grffile}
\usepackage{subcaption}
\usepackage{amsmath}
\usepackage{amsthm}
\usepackage{amssymb}
\usepackage{bm}
\usepackage{booktabs}

\usepackage{natbib}
\usepackage{soul}
\usepackage[framemethod=tikz]{mdframed}
\usepackage{outlines}
\usepackage{amsmath}
\usepackage{mathtools}
\usepackage{mathrsfs}
\usepackage{cancel}
\usepackage{blkarray, bigstrut}
\usepackage[utf8]{inputenc}
\usepackage{cancel}
\usepackage[ruled,vlined]{algorithm2e}
\usepackage{algorithm2e}
\usepackage{tikz}
\usepackage{float}
\usepackage{bibunits}

\usepackage[section]{placeins}

\theoremstyle{plain}
\newtheorem{thm}{Theorem}
\newtheorem{lem}{Lemma}
\newtheorem{prop}{Proposition}
\newtheorem{cor}{Corollary}
\newtheorem{thma}{Theorem}

\DeclarePairedDelimiter\floor{\lfloor}{\rfloor}

\makeatletter
\def\Vec{\mathop{\operator@font vec}\nolimits}
\makeatother

\newsavebox{\fminipagebox}
\NewDocumentEnvironment{fminipage}{m O{\fboxsep}}
{\par\kern#2\noindent\begin{lrbox}{\fminipagebox}
		\begin{minipage}{#1}\ignorespaces}
		{\end{minipage}\end{lrbox}%
	\makebox[#1]{%
		\kern\dimexpr-\fboxsep-\fboxrule\relax
		\fbox{\usebox{\fminipagebox}}%
		\kern\dimexpr-\fboxsep-\fboxrule\relax
	}\par\kern#2
}

\defaultbibliography{rbLit}
\defaultbibliographystyle{plainnat}

\newcommand{\Xs}{\mathscr{X}}

\makeatletter
\def\ps@pprintTitle{%
	\let\@oddhead\@empty
	\let\@evenhead\@empty
	\def\@oddfoot{}%
	\let\@evenfoot\@oddfoot}
\makeatother

\begin{document}

    \title{\bf Control variates and Rao-Blackwellization for deterministic sweep Markov chains}
    \author{Stephen Berg, Jun Zhu, and Murray K. Clayton
      \hspace{.2cm}\\
      Department of Statistics, University of Wisconsin-Madison}
    \maketitle

	\begin{abstract}
		We study control variate methods for Markov chain Monte Carlo (MCMC) in the setting of deterministic sweep sampling using $K\geq 2$ transition kernels. New variance reduction results are provided for MCMC averages based on sweeps over general transition kernels, leading to a particularly simple control variate estimator in the setting of deterministic sweep Gibbs sampling. Theoretical comparisons of our proposed control variate estimators with existing literature are made, and a simulation study is performed to examine the amount of variance reduction in some example cases. We also relate control variate approaches to approaches based on conditioning (or Rao-Blackwellization), and show that the latter can be viewed as an approximation of the former. Our theoretical results hold for Markov chains under standard geometric drift assumptions. 
	\end{abstract}
    
\begin{bibunit}

\section{Introduction\label{sec:intro}}
	
Markov chain Monte Carlo (MCMC) is a widely used technique for drawing samples from intractable probability distributions. In statistics, MCMC is now a standard tool in Bayesian analysis for sampling from complicated posterior distributions. The goal of MCMC is usually to approximate quantities such as $\int \pi(dx)g(x)$, where $\pi$ is an intractable probability measure, and $g:X\to\mathbb{R}^d$ is a $\pi$-integrable function mapping a state space $X$ to $\mathbb{R}^d$ for some integer $d\geq 1$. In MCMC, a Markov chain $X_0,X_1,X_2,...$ with a stationary probability measure $\pi$ is simulated for some finite number of iterations $M$, and $\int \pi(dx)g(x)$ is then approximated by the empirical average \begin{align}
	S_M/M=M^{-1}\sum_{t=0}^{M-1}g(X_t).\label{eq:ordinary}
	\end{align} Since MCMC sampling is stochastic, the empirical average $S_M/M$ based on MCMC simulations is inherently variable. For finite $M$, $S_M/M$ will even be biased, except for the usually unrealistic case when the Markov chain is started at a random draw from the stationary distribution $\pi$. However, under suitable conditions, a central limit theorem can be shown for $S_M/M$ stating that \begin{align}
	M^{1/2}\left\{S_M/M-\int\pi(dx)g(x)\right\}\overset{d}{\to}N(0,\Sigma)\label{eq:clt}
	\end{align} as $M\to\infty$, where $\overset{d}{\to}$ denotes convergence in distribution~\citep{meyn2009markov}. In this sense, $S_M/M$ is asymptotically unbiased as $M\to\infty$, and the MCMC error asymptotically comes entirely from the variance $\Sigma$. Due to this fact, one sensible measure of the efficiency of a Markov chain averaging scheme is the asymptotic variance $\Sigma$, and we will use this notion of efficiency in the remainder.

A variety of techniques exist for reducing the variance $\Sigma$ in~\eqref{eq:clt} for MCMC simulations, including conditioning, control variates, and antithetic sampling \citep[see, e.g.,][]{liu2008monte,robert2004monte}. We focus on control variate approaches here, although we also make connections to conditioning based approaches. In control variate approaches, mean zero random variables are added to each term of \eqref{eq:ordinary} in such a way that the variance of the sum is reduced. In approaches based on conditioning, the overarching idea is to replace $g$ in~\eqref{eq:ordinary} with a conditional expectation $E(g|\mathscr{H})$, where $\mathscr{H}$ is a $\sigma$-field. The hope is that the resulting average may have a reduced variance relative to~\eqref{eq:ordinary} due to Jensen's inequality, similar to the classical Rao-Blackwell theorem~\citep{rao1945information,blackwell1947}. This similarity has led to the common use of the term Rao-Blackwellization to describe techniques in which a MCMC average of a conditional expectation is taken in order to reduce the asymptotic variance $\Sigma$. However, in MCMC, unlike in classical Monte Carlo settings, independence does not hold and a naive conditioning approach may increase the asymptotic variance \citep{geyerConditioning}.

Previous results on variance reduction for MCMC apply primarily to Markov chain samplers with time-homogeneous transitions. Many of these results apply specifically to reversible Markov chains, such as the ones resulting from Metropolis samplers, Metropolis-Hastings samplers, or random sweep versions of Gibbs or Metropolis samplers~\citep[][]{metropolis1953equation,hastings1970monte,geman1984stochastic}. For example, \citet{casella1996rao} provided variance reduction results for Markov chains resulting from the Metropolis-Hastings algorithm. A conditioning approach was used in \citet{mckeague} to obtain a variance reduction result for reversible Markov chains. In \citet{meyn2008control}, control variate methods were discussed for time-homogeneous Markov chains in the context of network models. \citet{douc2011} studied a Rao-Blackwellization method for Markov chains based on Metropolis-Hastings algorithms. In \citet{dellaportas2012control}, a control variate method was given for reducing the variance of estimates based on reversible Markov chains. In \citet{brosse2018diffusion}, a control variate scheme was used to obtain variance reductions for certain Markov chains that can be related through a limiting process to a Langevin diffusion.

In this work, we study control variate methods for deterministic sweep sampling schemes, where multiple Markov transition kernels are cycled through in a fixed order. Often, deterministic sweep schemes arise when a vector valued Markov chain state is updated component-wise, and a single sweep updates every component. Deterministic sweep schemes include, but are not limited to, the commonly used deterministic sweep Gibbs sampler. 

The main technical difference between the deterministic sweep setting here and the setting of previous work is that deterministic sweep samplers are (for $K\geq 2$ components) time-inhomogeneous and not reversible, even when the individual updates are themselves reversible. To our knowledge, the major previous deterministic sweep variance reduction result comes from \citet{wongKongLiu}. Their result applies in the special setting of a Gibbs sampler on a state space with two components, where the integrand $g$ being averaged only depends on one of the components. In contrast, all of our results apply to integrands $g$ with arbitrary functional dependence on the components of the state. Furthermore, we show that the conditioning estimate proposed in \citet{wongKongLiu} can be viewed as an instance of our control variate scheme with a suboptimal choice of weight, and that choosing the weight optimally leads to a lower asymptotic variance. Unlike the estimator in \citet{wongKongLiu}, our proposed control variate estimator uses the values of both components in the Gibbs sampler chain.

Our methodology adapts the control variate basis function approach of \citet{dellaportas2012control} to our setting of time-inhomogeneous, deterministic sweep Markov chains, from the original setting of time-homogeneous, reversible Markov chains. While this change in setting requires technical modifications, it also leads to a useful practical consequence: any two-component deterministic sweep Gibbs sampler with the control variate scheme described here is guaranteed to achieve a lower variance than the analogous two-component random sweep Gibbs sampling approach from \citet{dellaportas2012control}. Two-component Gibbs samplers have long been known to have lower asymptotic variance for deterministic sweep than for random sweep \citep{greenwood1998}, and our result extends this earlier work to the control variate setting, where (for deterministic sweep) the integrand is time-inhomogeneous.

In Section~\ref{sec:mcmcBasics}, we introduce more details about MCMC. In Section~\ref{sec:theoryResults}, we introduce our technical assumptions and state our variance reduction results. In Section~\ref{sec:comparisons}, we provide theoretical comparisons between the control variate methodology here and estimators from prior literature. In Section~\ref{sec:numerical}, we describe some example cases to which our methodology applies and perform a simulation study. We conclude with a discussion of applications and potential extensions in Section~\ref{sec:conclusion}.

	\section{MCMC basics, motivation, and notation\label{sec:mcmcBasics}}


	\subsection{$K$-component samplers}
	For integers $K\geq 1$, we define $K$-component deterministic and random sweep samplers. We consider Markov chains evolving on a state space $(X,\mathscr{X})$, where $X$ is assumed to be a complete separable metric space, and $\mathscr{X}$ is the associated Borel $\sigma$-algebra. We say a function $\Pi:X\times\mathscr{X}\to [0,1]$ is a probability kernel if $\Pi(\cdot,A):X\to [0,1]$ is an $\mathscr{X}$--measurable function of $x$ for each $A\in\mathscr{X}$, and also $\Pi(x,\cdot):\mathscr{X}\to [0,1]$ defines a probability measure on $(X,\mathscr{X})$ for each $x\in X$. Given a probability measure $\lambda$ on $(X,\mathscr{X})$, we say a probability kernel $\Pi(x,A):X\times \mathscr{X}\to [0,1]$ is $\lambda$-stationary iff $\lambda(A)=\int\lambda(dx)\Pi(x,A)$ for all $A\in\mathscr{X}$. 
	
	We use $\mathbb{R}$ to denote the extended real line $\left[-\infty,\infty\right]$. For functions $f:X\to\mathbb{R}^p$ where $p\geq 1$ and probability kernels $\Pi:X\times \mathscr{X}\to [0,1]$, we define $\Pi^0f(x)=f(x)$, $\Pi^1f(x)=\Pi f(x)=\int\Pi(x,dy) f(y)$, and $\Pi^tf(x)=\Pi(\Pi^{t-1}f)(x)$ for $t>1$. We also define $\Pi^t(x,A)$ by $\Pi^t(x,A)=\Pi^tI_A(x)$, where $I_A(\cdot):X\to\mathbb{R}$ denotes the indicator function with $I_A(x)=1$ for $x\in A$ and $I_A(x)=0$ elsewhere. 

	Further, we define the sequence space $(\Omega,\mathscr{F})$ by $\Omega=X^{\mathbb{N}}$ and the product $\sigma$-algebra $\mathscr{F}$. We define a Markov chain as a collection $\{X_t\}_{t=0}^{\infty}$ of random variables $X_t:\Omega\to X,\omega\to\omega_t$ for $\omega\in \Omega$ with joint law $P_\nu$ specified by an initial probability measure $\nu$ on $(X,\mathscr{X})$ and a set of probability kernels $\{\Pi_t:X\times \mathscr{X}\to [0,1]\}_{t=1}^{\infty}$ through \begin{align}
		&P(X_0\in A_0,X_1\in A_1,...,X_t\in A_t)\\&=\int \nu(dx_0)\Pi_1(x_0,dx_1)...\Pi_t(x_{t-1},dx_t)I(\{x_i\in A_i,\forall i\})\nonumber
	\end{align} for any finite collection $A_0,A_1,...,A_t\in \mathscr{X}$. Under our topological assumptions on $(X,\mathscr{X})$, the Kolmogorov extension theorem guarantees that such a law $P_\nu$ exists and is unique \citep{durrett2010probability}.

	We consider deterministic sweep and random sweep samplers induced by a family $\Pi_k:\;X\times \Xs\to[0,1],k=1,...,K$ of $\pi$-stationary Markov transition kernels. The random sweep sampling scheme updates the chain at each step $t$ via the Markov transition kernel $Q_t=Q=K^{-1}\sum_{k=1}^{K}\Pi_k$. 
	

	For deterministic sweep sampling schemes, we need some more notation. We define the permutation function $\sigma(\cdot): \{1,...,K\}\to\{1,...,K\}$ by $
	\sigma(k)=k+1$ for $k<K$ and $\sigma(k)=1$ for $k=K$. We also define $\sigma^t(k)$ inductively by $\sigma^t(k)=k$ for $t=0$, and $\sigma^t(k)=\sigma\{\sigma^{t-1}(k)\}$ for $t>0$.

	Deterministic sweep samplers update the Markov chain by applying the kernels $\Pi_k$ in a fixed order. At time $0$, the transition operator used to update the state is, without loss of generality, $\Pi_1=\Pi_{\sigma^0(1)}$, and at time $t$, the transition operator is $\Pi_{\sigma^t(1)}$. Thus, for an initial probability measure $\nu$ on $(X,\mathscr{X})$, we have \begin{align*}
		&P_\nu(X_0\in A_0,...,X_t\in A_t)\\&=\int \nu(dx_0)\Pi_1(x_0,dx_1)...\Pi_{\sigma^t(1)}(x_{t-1},dx_t)I(\{x_i\in A_i,\forall i\}).
	\end{align*} We also define $P_k^{t}f(x)=\left\{\Pi_{\sigma^0(k)}\Pi_{\sigma^1(k)}\cdots \Pi_{\sigma^{t-1}(k)}f\right\}(x)$ and $P_k^t(x,A)=P_k^tI_A(x)$. In these notations, $P^t_k$ is a composition of multiple kernels rather than a repeated composition of a single kernel. 

	We will sometimes use the notation $\braket{f,g}$ in order to denote the integral $\int\pi(dx)f(x)g(x)$.


\subsection{Gibbs kernels\label{sec:gibbs}}

Now, we define the probability kernels used in Gibbs sampling. Define the identity map $Y:(X,\Xs)\to (X,\Xs), x\to x$ for $x\in X$. We first define a regular conditional distribution kernel. We say a probability kernel $\Pi:X\times \Xs\to [0,1]$ is a regular conditional distribution kernel with respect to $(\mathscr{G},\pi)$, where $\mathscr{G}\subset\mathscr{X}$ is a sub-$\sigma$-algebra of $\Xs$, whenever
\begin{itemize}
\item $\Pi(x,A)=E\{I(Y\in A)|\mathscr{G}\}$ almost everywhere with respect to $\pi$ (a.e. $\pi$), for each $A\in\mathscr{X}$ 
\item For $\pi$-a.e. $x$, $\Pi(x,A)$ is a probability measure on $(X,\mathscr{X})$
\end{itemize}

It is well-known that when $X$ is a complete separable metric space and $\mathscr{X}$ the associated Borel $\sigma$-algebra, such regular conditional distributions $\Pi$ always exist \citep[e.g.,][]{durrett2010probability}. 

We now define Gibbs kernels. Gibbs sampling is typically described and implemented in terms of conditional distributions with respect to the $\sigma$-algebras generated by random variables, rather than general $\sigma$-algebras. Given a measurable function $h:X\to\mathbb{R}^n$ for some $n\geq 1$, we say a probability kernel $\Pi:X\times \mathscr{X}\to [0,1]$ is a Gibbs kernel with respect to $(h,\pi)$ if $\Pi$ is a regular conditional distribution kernel with respect to the $\sigma$-algebra $\sigma(h)$.



We now provide a concrete example for a bivariate normal stationary distribution $\pi$. Let $X=\mathbb{R}^2$ and let $\mathscr{R}^2$ be the Borel $\sigma$-algebra on $\mathbb{R}^2$. Additionally, use $x=[x_1,x_2]^T$ to denote elements $x\in\mathbb{R}^2$, and let $\pi$ refer to the probability measure on $(X,\mathscr{X})$ such that the identity map $Y:(\mathbb{R}^2,\mathscr{R}^2)\to (\mathbb{R}^2,\mathscr{R}^2),x\to x$ follows a normal distribution with mean $\mu=[0,0]^T$ and covariance $\Sigma=\begin{bmatrix}1 & \rho \\ \rho & 1\end{bmatrix}$. Define the coordinate projections $h_i:(\mathbb{R}^2,\mathscr{R}^2)\to (\mathbb{R},\mathscr{R}),[x_1,x_2]^T\to x_i$ for $i=1,2$. Use $\delta_{x_1}$ and $\pi_{x_1}$ to denote the probability measures on $(\mathbb{R},\mathscr{R})$ with unit mass at $x_1$, and for a normal random variable with mean $\rho x_1$ and variance $1-\rho^2$, respectively. Then the probability kernels $\Pi_i(x,A):\mathbb{R}^2\times \mathscr{R}^2\to [0,1]$ defined by
	\begin{align}
		&\Pi_1(x,A)=\int (\delta_{x_1}\times \pi_{x_1})(dx'_1,dx'_2)I(x'\in A)\label{eq:pi1}\\ 
		&\Pi_2(x,A)=\int (\pi_{x_2}\times \delta_{x_2})(dx'_1,dx'_2)I(x'\in A)\label{eq:pi2}
	\end{align}
are Gibbs kernels with respect to $(h_i,\pi)$.\qed

Gibbs kernels have some useful properties. First, any Gibbs kernel $\Pi$ with respect to $(h,\pi)$ preserves $\pi$. This follows since \begin{align}
	\int\pi(dx)\Pi(x,A)=\pi(A)
	\end{align} for each $A\in\mathscr{X}$, from the properties of conditional expectation. Additionally, the idempotence property $\Pi\{\Pi f\}(x)=\Pi f(x)$ holds a.e. $\pi$ for each $\pi$-integrable $f$. Finally, for functions $f,g$ which are square-integrable with respect to $\pi$, we have \begin{align}
	&\int \pi(dx)f(x)\Pi g(x)=\int \pi(dx)\Pi \{f\Pi g\}(x)\nonumber\\&=\int \pi(dx)\Pi f(x)\Pi g(x)=\int\pi(dx)\Pi f(x)g(x).\label{eq:reversible}
	\end{align} The equality $\braket{f,\Pi g}=\braket{\Pi f,g}$ is the useful \textit{reversibility} property. Thus from~\eqref{eq:reversible}, we see that Gibbs transition kernels are each reversible with respect to $\pi$. However, compositions of the form $\Pi_1\Pi_2\cdots \Pi_k$ will in general not be reversible with respect to $\pi$.

\subsection{Control variates and Rao-Blackwellization}

We will henceforth use $g:X\to\mathbb{R}^d$ to denote a function for which we would like to obtain $\int\pi(dx)g$ using MCMC, and $f:X\to\mathbb{R}^p$ to denote a control variate basis function \citep{dellaportas2012control}. For deterministic sweep chains $\{X_t\}_{t=0}^{\infty}$ that use probability kernel $\Pi_{\sigma^{t}(1)}$ to obtain $X_t$, we consider four averaging schemes:
\begin{flalign}
	&\text{Empirical: }&\label{eq:emp}\\&M^{-1}S_M=M^{-1}\sum_{t=0}^{M-1}g(X_t)&\nonumber
\end{flalign}
\begin{flalign}
	&\text{Rao-Blackwellized: }&\label{eq:rb}\\&M^{-1}S_M=M^{-1}\sum_{t=0}^{M-1}\Pi_{\sigma^t(1)}g(X_t)&\nonumber
\end{flalign}
\begin{flalign}
	&\text{Fixed weight control variate: } &\label{eq:cv0}\\&M^{-1}S_M=M^{-1}\sum_{t=0}^{M-1}g(X_t)-C^T\{f(X_t)-\Pi_{\sigma^t(1)}f(X_t)\}\nonumber
\end{flalign}
\begin{flalign}
	&\text{General control variate: }\label{eq:cv1}&\\
	&M^{-1}S_M=\sum_{t=0}^{M-1}g(X_t)-C_{\sigma^t(1)}^Tf(X_t)+C_{\sigma^{t+1}(1)}^T\Pi_{\sigma^t(1)}f(X_t)\nonumber
\end{flalign} In~\eqref{eq:cv0} and~\eqref{eq:cv1}, we use $C$ and $C_k$ for $\;k=1,...,K$ to refer to fixed $p\times d$ weight matrices. We will show for suitably chosen $C$ ($C_k$), the control variate approaches have a smaller asymptotic variance than the empirical average. We note that the Rao-Blackwellized estimator \eqref{eq:rb} is a special case of the control variate estimators, using the choices $f=g$ and $C_k=I_{d\times d}$. 



Because the control variate and Rao-Blackwellization schemes involve the conditional expectations $\Pi_kf$ (or $\Pi_kg$, in the Rao-Blackwellization setting), it is necessary in practice for these conditional expectations to have a computationally tractable form in order to apply the control variate and Rao-Blackwellization schemes.

\section{Assumptions and variance reduction results\label{sec:theoryResults}}

In this section, we consider Markov chains $\{X_t\}_{t=0}^{\infty}$ with arbitrary initial law $\nu$ on $(X,\Xs)$ and time-inhomogeneous transition kernels $P_t=\Pi_{\sigma^t(1)}$. Throughout, we will take $\pi$ to be a probability measure on $(X,\Xs)$ which we would like to take expectations with respect to. 

	
\subsection{Assumptions\label{subSec:assumptions}}

First, we define small sets, in the sense of \citet{meyn2009markov}. We say a set $C\in\mathscr{X}$ is a small set with respect to a probability kernel $P$ if there exists an $m>0$, and a non-trivial measure $\nu_m$ on $(X,\mathscr{X})$, such that \begin{align*}
	P^m(x,B)\geq \nu_m(B)
\end{align*} for all $x\in C$, $B\in\mathscr{X}$.

We make the following assumptions on the composition kernels $P_k^{Kt}$ and the transition kernels $\Pi_1,...,\Pi_K$.

\begin{enumerate}[label=(A.\arabic*)]
	\item \label{cond:piStationary} The kernels $\Pi_k$ are $\pi$-stationary. This holds whenever at least one of \ref{cond:piReversible} and \ref{cond:projection} holds:
	\begin{enumerate}[label=(A.\arabic{enumi}-\arabic*)]
	\item \label{cond:piReversible} (Reversibility) The $\Pi_k$ are reversible, so that $\braket{f,\Pi_kg}=\braket{\Pi_kf,g}$ for all square integrable functions $f,g:X\to\mathbb{R}$.
	\item \label{cond:projection} (Gibbs kernels) The $\Pi_k$ are Gibbs kernels with respect to $(h_k,\pi)$ for some set of functions $h_k:(X,\Xs)\to(\mathbb{R}^{m_k},\mathscr{R}^{m_k})$, where $m_k\geq 1$ is an integer denoting the dimension of the range of $h_k$.
	\end{enumerate}
	\item \label{cond:irreducibility}($\psi$-irreducibility) There exists a probability measure $\psi$ on $(X,\mathscr{X})$ such that for each $k=1,...,K$ and all $A\in \mathscr{X}$ with $\psi(A)>0$, and for all $x\in X$, there exists a positive integer $t=t(x,A,k)$ such that $P_k^{Kt}(x,A)>0$.
	\item \label{cond:drift}(Geometric drift) There exist small sets $C_k\in\mathscr{X}$ with respect to $P_k^K$, constants $\lambda_k<1$ and $b_k<\infty$, and  functions $V_k:\;X\to [1,\infty)$, such that for $k=1,...,K$, \begin{align*}
		P_k^K V_k(x)\leq \lambda_kV_k(x)+b_kI_{C_k}(x).
	\end{align*}
	\item \label{cond:aperiodicity} (Aperiodicity) The composition kernels $P_k^K$ are assumed to be aperiodic.
\end{enumerate}

Assumption~\ref{cond:piStationary} ensures that the transition kernels $\Pi_k$ preserve the stationary distribution $\pi$. Assumption~\ref{cond:irreducibility} ensures that $\pi$ is the unique stationary distribution for the $K$-step composition kernels $P_k^K$ for $k=1,...,K$ \citep[see, e.g.][]{meyn2009markov}. Assumption~\ref{cond:drift} 
ensures that Markov chains with transition kernel $P_k^K$ are Harris recurrent, so that for any $A\in\mathscr{X}$ with $\psi(A)>0$, we have $P_x(\cap_{N=1}^{\infty}\cup_{k=N}^{\infty}\{X_t\in A\})=1$ for all $x
\in A$, where $P_x$ refers to the Markov chain law with point mass initial distribution $\delta_x$ and transition kernel $P_k^{K}$. Regarding~\ref{cond:aperiodicity}, we say an irreducible probability kernel is aperiodic if $d=1$ is the largest integer such that there exist sets $D_1,...,D_d$ satisfying \begin{enumerate}
	\item $P_k^K(x,D_{i+1})=1$ for $x\in D_i$, $i\equiv 0,...,d-1$ (mod $d$)
	\item $\psi\{(\cup_{i=1}^{d}D_i)^C\}=0$
	\item $D_1,...,D_{d}$ are disjoint
\end{enumerate} \citep[see, e.g.][]{meyn2009markov}. Lemma~\ref{lem:aperiodic} in Appendix~\ref{app:appA} shows that Assumption~\ref{cond:projection}, when it holds, ensures~\ref{cond:aperiodicity} holds also.

We make the following assumptions about the functions $g:X\to\mathbb{R}^d$, $f:X\to\mathbb{R}^p$, and $V_k(x)$ in~\ref{cond:drift}:

\begin{enumerate}[label=(B.\arabic*)]
	\item \label{cond:integrability} (square integrability)$\int \pi(dx)V_k^2(x)<\infty$
	\item \label{cond:smallg} $|a^Tg(x)|\leq V_k(x)$ for all $a\in\mathbb{R}^d$ with $\|a\|_2\leq 1$, where $\|x\|_2$ denotes the Euclidean norm of $x$. For a univariate function $g$, this is equivalent to assuming that $|g(x)|\leq V_k(x)$ for all $x$.
	\item \label{cond:centered}$\int\pi(dx)g(x)=0$
	\item \label{cond:smallf}$\int \pi(dx)f^Tf<\infty$ and $f^Tf<\infty$ for all $x\in X$
\end{enumerate}

Assumptions~\ref{cond:drift},~\ref{cond:aperiodicity},~\ref{cond:integrability}, and~\ref{cond:smallg} will be used to ensure certain bounds on solutions $\tilde{g}_k$ to the Poisson equations $\tilde{g}_k-P_k^{K}\tilde{g}_k=g-\int\pi(dx)g(x)$. Assumption~\ref{cond:centered} is introduced for notational convenience: it allows us to write, for example, statements such as $\tilde{g}_k-P_k^{K}\tilde{g}_k=g$ rather than $\tilde{g}_k-P_k^{K}\tilde{g}_k=g-\int\pi(dx)g(x)$. For a general integrand $g$, the results to follow will apply to the function $\bar{g}(x)=g(x)-\int\pi(dx')g(x')$. We note that \ref{cond:smallf} will be satisfied for $f=g$ when~\ref{cond:piStationary}--\ref{cond:aperiodicity} and~\ref{cond:integrability}--\ref{cond:smallg} hold.

\subsection{Results\label{subSec:main1}}

In this section, we state our main results, and defer the proofs to the appendices.

\begin{prop}\label{prop:poisEq}
Under~\ref{cond:piStationary}--\ref{cond:aperiodicity} and~\ref{cond:integrability}--\ref{cond:centered}, there exist functions $\hat{g}_k:X\to\mathbb{R}^d$ with \begin{align}
&\hat{g}_k(x)=\sum_{t=0}^{\infty}P_{k}^tg(x)\;\;\;\;\;\;k=1,...,K\label{eq:gdef1}
\end{align} a.e. $\pi$. The sums in the definition of $\hat{g}_k$ are absolutely convergent elementwise for $\pi$-a.e. $x\in X$, and $\int \pi(dx)\hat{g}_k^T\hat{g}_k<\infty$ for each $k$. Additionally, each $\hat{g}_k$ satisfies a corresponding Poisson-type equation \begin{align*}
	\hat{g}_k-\Pi_k\hat{g}_{\sigma(k)}=g,\;\;\;a.e.\;\pi.
\end{align*}
\end{prop}

Proposition~\ref{prop:poisEq} is proved in Appendix~\ref{app:appB}.

We will use the Poisson equation solutions from Proposition~\ref{prop:poisEq} in order to write each of the sums $S_M$ in~\eqref{eq:emp}-\eqref{eq:cv1} as the sum of an approximating martingale, plus a small error term. We will then apply central limit theorems for martingales to obtain results for the asymptotic variance of MCMC averages~\citep{gordin1969central,hall1980martingale}. The Poisson equation technique is discussed for time-homogeneous Markov chains in Chapter 17 of \citet[][]{meyn2009markov}. We also found useful ideas in \citet{Ynes1988}, in a stochastic approximation setting involving a deterministic sweep Gibbs sampler. We modify these approaches to fit the setting here, where both the transition kernels $\Pi_k$, and the integrand, are time inhomogeneous.

Intuitively, one can verify that $\hat{g}_k-\Pi_k\hat{g}_{\sigma(k)}=g$ by checking that each term in $\hat{g}_k$ matches a term in $\Pi_k\hat{g}_{\sigma(k)}$ except for the first term $g$, so that all terms besides $g$ cancel, provided the sums and integrals can be rearranged as needed.
 We define \begin{align}
	&U_k=\int\pi(dx)\{ff^T-(\Pi_kf)(\Pi_kf)^T\}\label{eq:Uk}\\
	&V_k=\int \pi(dx)\{f\hat{g}_{\sigma(k)}^T-(\Pi_kf)(\Pi_k\hat{g}_{\sigma(k)}^T)\}.\label{eq:Vk}
\end{align} 

We now consider sums of the form \begin{align}
S_M=\sum_{t=0}^{M-1}g(X_t)-C_{\sigma^t(1)}^Tf(X_t)+C_{\sigma^{t+1}(1)}^T\Pi_{\sigma(t)}f(X_t),\label{eq:vecS}\end{align} where $C_k$ are $p\times d$ matrices of weights, for $k=1,...,K$. We establish the following result.
\begin{thm}\label{thm:vecG}
	Assume~\ref{cond:piStationary}--\ref{cond:aperiodicity} and~\ref{cond:integrability}--\ref{cond:smallf}. 
We have \begin{align*}M^{-1/2}S_M\overset{d}{\to}Z\end{align*} where $Z$ is a random $d$-vector with characteristic function $\exp\{-t^T\Sigma_C t/2\}$. The variance $\Sigma_C$ can be written as \begin{align}
		\Sigma_C&=\int \pi(dx)gg^T+K^{-1}\sum_{k=1}^{K}\sum_{t=1}^{\infty}\int\pi(dx)\{g(P_kg)^T+(P_kg)g^T\}\label{eq:SigmaC}\\
		&+K^{-1}\sum_{k=1}^{K}C_{\sigma(k)}^TU_kC_{\sigma(k)}-C_{\sigma(k)}^TV_k-V_k^TC_{\sigma(k)}.\nonumber
	\end{align} The variance $\Sigma_C$ is minimized\footnote{For positive semidefinite matrices $A$ and $B$, we say $A\geq B$ if $A-B$ is positive semidefinite. We say $A>B$ if $A-B$ is positive semidefinite with at least 1 nonzero eigenvalue.} at $C_{\sigma(k)}=\tilde{C}_{\sigma(k)}$, where $\tilde{C}_{\sigma(k)}=U_k^{\dagger}V_k$ and $A^{\dagger}$ denotes the pseudoinverse of $A$.
\end{thm}

The proof of Theorem~\ref{thm:vecG} is in Appendix~\ref{app:appB}.

Theorem~\ref{thm:vecG} states that $M^{-1/2}S_M\overset{d}{\to}N(0,\Sigma_C)$ whenever $\Sigma_C$ is positive definite. In general, the optimal weight expression $\tilde{C}_{\sigma(k)}=U_k^{\dagger}V_k$ in Theorem~\ref{thm:vecG} appears daunting, since the $V_k$ contain integrals involving the Poisson equation solutions $\hat{g}_{\sigma(k)}$. However, Corollary~\ref{cor:aVar0} below will show that the expression for the optimal weight is much simpler in the setting of deterministic sweep Gibbs sampling with fixed control variate weights $C_1=\cdots=C_K=C$. In particular, the optimal weight expression no longer explicitly involves the Poisson equation solutions.

We define \begin{align}
	U=K^{-1}\sum_{k=1}^{K}U_k\label{eq:U}\\
	V=K^{-1}\sum_{k=1}^{K}V_k\label{eq:V}
\end{align} where $U_k$ and $V_k$ are defined in~\eqref{eq:Uk} and~\eqref{eq:Vk}, respectively, for $k=1,...,K$.

We then have Corollaries~\ref{cor:aVar0} and~\ref{cor:naiveVec}.

\begin{cor}\label{cor:aVar0}
Suppose \ref{cond:piStationary},~\ref{cond:irreducibility}--\ref{cond:drift}, and~\ref{cond:integrability}--\ref{cond:smallf} hold. Then for the fixed weight scheme with $C_1=\cdots=C_K=C$, we have \begin{align*}
	\Sigma_C&=\int \pi(dx)gg^T+K^{-1}\sum_{k=1}^{K}\sum_{t=1}^{\infty}\int\pi(dx)\{g(P_kg)^T+(P_kg)g^T\}\\
	&+C^TUC-C^TV-V^TC
\end{align*}
	and $\Sigma_{C}$ is minimized at $C=\tilde{C}$, where $\tilde{C}=U^{\dagger}V$. If~\ref{cond:projection} also holds, then we have the simplified expression $V=\int\pi(dx)fg^T$ for $V$.
\end{cor}

Proof: A detailed proof is in Appendix~\ref{app:appB}. We outline the steps to obtain the simplified representation for $V$ under~\ref{cond:projection} here. We first observe $\int \pi(dx)(\Pi_kf)(\Pi_k\hat{g}_{\sigma(k)}^T)=\int \pi(dx)f\Pi_k\hat{g}_{\sigma(k)}^T$ by the idempotence and reversibility of $\Pi_k$ under~\ref{cond:projection}. Then, we use Proposition~\ref{prop:poisEq} and rearrange a sum to obtain \begin{align*}
	V&=K^{-1}\sum_{k=1}^{K}\int\pi(dx)\{f\hat{g}_{\sigma(k)}-f\Pi_k\hat{g}_{\sigma(k)}\}\\
	&=K^{-1}\sum_{k=1}^{K}\int\pi(dx)f(\hat{g}_{k}-\Pi_k\hat{g}_{\sigma(k)})^T=\int\pi(dx)fg^T.
\end{align*}\qed

Corollary~\ref{cor:aVar0} shows that the formula for the optimal control variate weight simplifies dramatically in the Gibbs sampling setting where~\ref{cond:projection} holds. In general, the quantity $V=K^{-1}\sum_{k=1}^{K}f\hat{g}_{\sigma(k)}^T-(\Pi_kf)(\Pi_k\hat{g}_{\sigma(k)})^T$ depends on the Poisson equation solutions $\hat{g}_k$ from Proposition~\ref{prop:poisEq}, whereas under~\ref{cond:projection}, the formula for $\tilde{C}$ no longer explicitly involves the $\hat{g}_{\sigma(k)}$. The proof obtaining the simplified representation for $V$ exploits both the idempotence and the reversibility properties of the $\Pi_k$ under~\ref{cond:projection}, in order to cancel all but the leading terms in the Poisson equation solutions $\hat{g}_{\sigma(k)}$. A similar cancellation also occurs in the time-homogeneous, reversible Markov chain setting~\citep{dellaportas2012control}.

In contrast, for general transition kernels or Gibbs sampling without fixed control variate weights, the optimal weights depend on high-order autocovariances due to the presence of the Poisson equation solutions $\hat{g}_k$ in $V$ ($V_k$). In the general setting, we will investigate an approach based on batch means to estimate the optimal control variate weights \citep[see, e.g.,][]{meyn2008control,dellaportas2012control}. Finally, we note that the matrix $U$ will always be invertible in the ordinary sense except when $a^Tf=b$ a.e. $\pi$ for some constant $b\in\mathbb{R}$ and nonzero $a\in\mathbb{R}^d$. This follows from Lemma~\ref{lem:smallCondExp} in Appendix~\ref{app:appA}.

We now consider deterministic sweep Gibbs sampling with the basis function $f=g$, for the fixed weight average $S_M/M=M^{-1}\sum_{t=0}^{M-1}g(X_t)-C^T\{g(X_t)-\Pi_{\sigma^t(1)}g(X_t)\}$. Let $0_{d\times d}$ and $I_{d\times d}$ refer to the $d\times d$ matrices of all $0$'s and the identity matrix, respectively, and use $\Sigma_0$ and $\Sigma_1$ to denote the variance $\Sigma_C$ when $C=0_{d\times d}$ and $C=I_{d\times d}$. Taking $C=I_{d\times d}$ yields $S_M=\sum_{t=0}^{M-1}\Pi_{\sigma^t(1)}g(X_t)$, which is exactly the Rao-Blackwellized average \eqref{eq:rb}. On the other hand, taking $C=0_{d\times d}$ is equivalent to the ordinary empirical average~\eqref{eq:emp}.

With this setup, we have the following result.

\begin{cor}\label{cor:naiveVec}
	Suppose \ref{cond:projection},~\ref{cond:irreducibility}--\ref{cond:drift}, and~\ref{cond:integrability}--\ref{cond:centered} hold. We have \begin{align*}&\Sigma_0=K^{-1}\sum_{k=1}^{K}\left[\int \pi(dx)gg^T+\sum_{t=1}^{\infty}\int\pi(dx)\{g(P_k^tg)^T+(P_k^tg)g^T\}\right]\end{align*}and \begin{align*}&\Sigma_1=\Sigma_0-\int\pi(dx)gg^T-K^{-1}\sum_{k=1}^{K}\int\pi(dx)(\Pi_kg)(\Pi_kg)^T\leq \Sigma_0.\end{align*}
\end{cor}

Proof: A detailed proof is in Appendix~\ref{app:appB}. The result follows from collecting terms and simplifying the variance from Theorem~\ref{thm:vecG} with $C=I_{d\times d}$ and $f=g$. We again exploit~\ref{cond:projection} to use the simplified formula for $V$ from Corollary~\ref{cor:aVar0}. When $f=g$, then $V=\int\pi(dx)fg^T=\int\pi(dx)gg^T$ . \qed

Corollary~\ref{cor:naiveVec} shows that the asymptotic variance is always smaller for the Rao-Blackwellized average $\Sigma_1$ than for the empirical average $\Sigma_0$. Thus, for deterministic sweep Gibbs sampling, the apparently naive Rao-Blackwellization strategy of averaging the conditional expectation of the integrand $g$, with respect to whichever transition kernel is being used to update $X_t$, leads to an improved asymptotic variance.

\subsection{Estimating the optimal control variate weight}
For the control variate schemes, we will usually need to estimate $U_k$ and $V_k$, or $U$ and $V$, based on MCMC, and implement the control variate estimate as a post-processing of the empirical estimator. 
Define $\bar{g}_M(X_t)=g(X_t)-M^{-1}\sum_{t'=0}^{M-1}g(X_{t'})$. For simplicity, we give expressions for estimates for the case $M=NK$ where $N$ is an integer corresponding to the number of complete sweeps. 

For the fixed weight control variate approaches, one can use \begin{align}
\hat{C}_{M}=\hat{U}_{M}^{\dagger}\hat{V}_{M}\label{eq:cEst}\end{align} for suitable $\hat{U}_M$ and $\hat{V}_M$. One estimator $\hat{U}_M$ of $U$ is \begin{align*}
&\hat{U}_M=(M-1)^{-1}\sum_{t=0}^{M-2}\{f(X_{t+1})-\Pi_{\sigma^t(1)}f(X_{t})\}\{f(X_{t+1})-\Pi_{\sigma^t(1)}f(X_{t})\}^T.
\end{align*} For Gibbs samplers, one can estimate $V$ by
\begin{align} 
	&\hat{V}_M=M^{-1}\sum_{t=0}^{M-1}f(X_t)\bar{g}_{M}(X_{t})^T.\label{eq:gibbsV}
\end{align} In general settings, we propose a batch means approach, where autocovariances up to a maximum lag of $B$ are included for some integer $B$: \begin{align}
	\hat{V}_M&=M^{-1}\sum_{t=0}^{M-1}f(X_t)\sum_{t'=t}^{(t+B)\wedge (M-1)}\bar{g}_M(X_{t'})^T\label{eq:genV}\\
	&-M^{-1}\sum_{t=0}^{M-1}\Pi_{\sigma^t(1)}f(X_t)\sum_{t'=t+1}^{(t+1+B)\wedge (M-1)}\bar{g}_M(X_{t'})^T,\nonumber
\end{align}

For $U_k$, we propose the estimates
\begin{align*}
	\hat{U}_{k,M}=(N-1)^{-1}\sum_{n=0}^{N-2}&\{f(X_{k+Kn})-\Pi_kf(X_{k+Kn-1})\}\\&\times\{f(X_{k+Kn})-\Pi_kf(X_{k+Kn-1})\}^T.
\end{align*} For $V_k$, we propose a batch means approach:
\begin{align}
	\hat{V}_{k,M}&=N^{-1}\sum_{n=0}^{N-1} f(X_{k+Kn})\sum_{t=k+Kn}^{(k+Kn+B)\wedge (NK-1)}\bar{g}_M(X_{t})^T\label{eq:varyingV}\\
	&-N^{-1}\sum_{n=0}^{N-1} \Pi_kf(X_{k+Kn-1})\sum_{t=k+Kn+1}^{(k+Kn+1+B)\wedge (NK-1)}\bar{g}_M(X_{t})^T.\nonumber
\end{align} We then compute estimates $\hat{C}_{\sigma(k),M}=\hat{U}_{k,M}^{\dagger}\hat{V}_{k,M}$. 

\section{Theoretical comparisons\label{sec:comparisons}}
\subsection{Comparison to \citet{wongKongLiu}\label{subSec:lwk}}

We next compare the control variate approaches studied here to the Rao-Blackwellization estimator from \citet{wongKongLiu}. The estimator in \citet{wongKongLiu} applies in the so-called \textit{data augmentation} Gibbs sampling setting, in which a two-component Gibbs sampler has components $\Pi_1$ and $\Pi_2$, and the integrand $g$ satisfies $\Pi_kg=g$ a.e. $\pi$ for at least one of the components $k$.
We assume the integrand $g$ satisfies \begin{enumerate}[label=(C.\arabic*)]
	\item \label{cond:subspace} (data augmentation) $\Pi_2g=g$ a.e. $\pi$.
	\item \label{cond:variance} (non-degeneracy) $\int \pi(dx)gg^T$ is a positive definite matrix.
\end{enumerate} In this setting, the estimator from \citet{wongKongLiu} is \begin{align}
	S_M/M=M^{-1}\sum_{t=0}^{M-1}\Pi_1g(X_t).\label{eq:condEst}
  \end{align} 
  

  Denote by $\Sigma_0$, $\Sigma_1$, $\Sigma_2$, $\Sigma_{\text{LWK}}$, and $\Sigma_{\tilde{C}}$ the variances obtained by the empirical estimator~\eqref{eq:emp}, the Rao-Blackwellized estimator~\eqref{eq:rb}, the fixed weight control variate scheme with $C=2I_{d\times d}$, the conditioning estimate in~\eqref{eq:condEst} due to \citet{wongKongLiu}, and the fixed weight control variate scheme with the optimal weight $C=\tilde{C}$, respectively. Additionally, define $A=\int\pi(dx)gg^T$ and $B=\int\pi(dx)(\Pi_1g)(\Pi_1g)^T$. Note that $A$ and $B$ are symmetric. Then we have the following result.

  \begin{prop}\label{prop:lwkProp}
	Assume $K=2$, and that \ref{cond:projection}, \ref{cond:irreducibility}--\ref{cond:drift},~\ref{cond:integrability}--\ref{cond:centered}, and \ref{cond:subspace}--\ref{cond:variance} hold. Then $\Sigma_{2}=\Sigma_{\text{LWK}}$, and $\Sigma_{\tilde{C}}\leq \Sigma_{2}< \Sigma_1< \Sigma_0$, with \begin{align*}
		&\Sigma_{\tilde{C}}-\Sigma_2=-2B(A-B)^{-1}B\leq 0\\
		&\Sigma_2-\Sigma_1=-(A+3B)/2<0\\
		&\Sigma_1-\Sigma_0=-(B+3A)/2<0
	\end{align*}
  \end{prop}

  Proof: deferred to Appendix~\ref{app:appB}.

  In Proposition~\ref{prop:lwkProp}, an ordering is given for the variances obtained by the empirical estimator, the Rao-Blackwellization estimator~\eqref{eq:rb}, the LWK estimator, and the optimal control variate estimator, in the data augmented Gibbs sampling setting. From Proposition~\ref{prop:lwkProp}, the variance from the optimal fixed weight control variate approach will be no larger, and will typically be smaller, than the variance from the approach of \citet{wongKongLiu}, which in turn will be lower than the variances from both the empirical estimator as well as the Rao-Blackwellization approach in Corollary~\ref{cor:aVar0}.
  

  \subsection{Comparison to random sweep Gibbs sampling approaches\label{subSec:rev}}

  We now compare the asymptotic variances resulting from deterministic and random sweep Gibbs sampling control variate schemes, for $K=2$. For kernels $\Pi_1$ and $\Pi_2$, we define the random sweep kernel $Q=(\Pi_1+\Pi_2)/2$. Additionally, we define the function \begin{align}
	h=g-C^T(f-Qf)\label{eq:hFunction}.
  \end{align} Then we have the following result.
  
  \begin{thm}{\label{thm:detVrev}}
	Suppose~\ref{cond:projection}--\ref{cond:drift} and \ref{cond:integrability}--\ref{cond:smallf} hold, and that the number of components $K=2$. Let an arbitrary weight matrix $C$ be given. Then the asymptotic variances $\Sigma_C$ and $\Sigma_C^{rev}$ from the fixed-weight deterministic sweep control variate and random sweep schemes with weight matrix $C$ have the alternate representations \begin{align}
		&\Sigma_{C}=\int\pi(dx)hh^T+\sum_{t=1}^{\infty}\int \pi(dx)h(Q^th)^T\label{eq:sigmaC_h}\\
		&\Sigma_{C}^{rev}=\int\pi(dx)hh^T+2\sum_{t=1}^{\infty}\int \pi(dx)h(Q^th)^T.\label{eq:sigmaC_rev}
	\end{align} Additionally, writing $\tilde{C}$ for the optimal control variate weight from the deterministic sweep scheme, and defining 
	$\bar{h}=g-\bar{C}^T(f-Qf)$, 
	where $\bar{C}$ is the optimal weight for the random sweep scheme, we have \begin{align*}
		&\Sigma_{\tilde{C}}-\Sigma_{\bar{C}}^{rev}=-(\bar{C}-\tilde{C})^TU^{\dagger}(\bar{C}-\tilde{C})-\sum_{t=1}^{\infty}\int\pi(dx)\bar{h}(Q^t\bar{h})^T\leq 0.\\
	\end{align*}
  \end{thm}

  Proof: deferred to Appendix~\ref{app:appB}.

  Theorem~\ref{thm:detVrev} shows that when $K=2$, it is better, in terms of the asymptotic variance, to use the deterministic sweep rather than the random sweep scheme. A similar result to Theorem~\ref{thm:detVrev} is given in \citet{greenwood1998}, but the \citet{greenwood1998} result applies to integrands that do not depend on the time step, unlike the control variate estimates~\eqref{eq:cv0} and~\eqref{eq:cv1}. In \citet{greenwood1998}, it is shown that for a fixed integrand $g$, the asymptotic variance resulting from deterministic sweep is always smaller than the asymptotic variance resulting from random sweep, for two-component Gibbs samplers. 


  \subsection{Discussion\label{subsec:Disc}}

We now show that the Rao-Blackwellization estimator~\eqref{eq:rb} can be viewed as an approximate control variate scheme. Suppose the control variate basis function $f=g$, and assume for simplicity that the matrices $U_k$ in~\eqref{eq:Uk} are positive definite. Then under the assumptions of Theorem~\ref{thm:vecG}, the optimal control variate weights are $\tilde{C}_{\sigma(k)}=U_k^{-1}V_k$, where $V_k=\int \pi(dx)\{g\hat{g}_{\sigma(k)}^T-(\Pi_kg)(\Pi_k\hat{g}_{\sigma(k)})^T\}$. Now, consider using the approximation \begin{align}
	\hat{g}_{\sigma(k)}=g+\Pi_{\sigma^1(k)}g+\Pi_{\sigma^1(k)}\Pi_{\sigma^2(k)}g+...\approx g,\label{eq:approxCV}
\end{align} within $V_k$, where the infinite sum in the Poisson equation solution is truncated after a single term.
Then we are left with \begin{align}
	V_k\approx \int\pi(dx)\{gg^T-(\Pi_kg)(\Pi_kg)^T\}=U_k
.\label{eq:approxCV2}\end{align}
 Thus, under the one term approximation of $\hat{g}_{\sigma(k)}$, we obtain $\tilde{C}_{\sigma(k)}=I_{d\times d}$ for each $k$. Similarly, for the fixed weight scheme, we obtain $\tilde{C}=I_{d\times d}$. But these choices of weights in~\eqref{eq:cv0} and~\eqref{eq:cv1} both lead to the estimator \begin{align*}
	S_M/M=M^{-1}\sum_{t=0}^{M-1}\Pi_{\sigma^t(1)}g(X_t),
\end{align*} which is exactly the Rao-Blackwellized estimator in~\eqref{eq:rb}.

An inspection of the expressions for the batch means estimates~\eqref{eq:genV} and~\eqref{eq:varyingV} shows that setting $B=0$ in the batch means estimate with $f=g$ can be interpreted as invoking the approximation~\eqref{eq:approxCV2}. In particular, in the setting where $f=g$ and $B=0$, we expect $\hat{U}_k\overset{a.s.}{\to}U_k$ and $\hat{V}_k\overset{a.s.}{\to}U_k$ asymptotically, so that the estimators from the control variate approach with $B=0$ will be identical to the estimates from the Rao-Blackwellized estimate~\eqref{eq:rb}. Figure~\ref{fig:isingBatch} in Section~\ref{sec:numerical} provides a numerical demonstration of this fact.

\section{Numerical examples\label{sec:numerical}}

\subsection{Bivariate normal}

First, we consider the bivariate normal Gibbs sampling setting from Section~\ref{sec:gibbs}. The deterministic sweep scheme updates $X_t$ by drawing from $\Pi_{1}(X_t,\cdot)$ for even $t$ and by drawing from $\Pi_{2}(X_t,\cdot)$ for odd $t$, where $\Pi_1$ and $\Pi_2$ are defined in \eqref{eq:pi1} and \eqref{eq:pi2}. It can be shown that for $k=1,2$, the composition kernels $P_k^2$ satisfy Assumption~\ref{cond:irreducibility} with the measure $\psi=\pi$. Additionally, Lemma~\ref{lem:aperiodic}, combined with~\ref{cond:projection} and~\ref{cond:irreducibility}, shows that~\ref{cond:aperiodicity} also holds. We verify in Lemma~\ref{lem:exampleDrift} (Appendix~\ref{app:appA}) that Assumption~\ref{cond:drift} also holds for $P_1^2$ and $P_2^2$ with the functions $V_1(x)=x_1^2+rx_2^2+1$ and $V_2(x)=rx_1^2+x_2^2+1$ for appropriately chosen $r>0$. These $V_k$ satisfy~\ref{cond:integrability}. Thus, Theorem~\ref{thm:vecG} can be applied to each of the following examples.

\begin{figure}[H]
\begin{subfigure}{0.3\textwidth}
\includegraphics[scale=0.22]{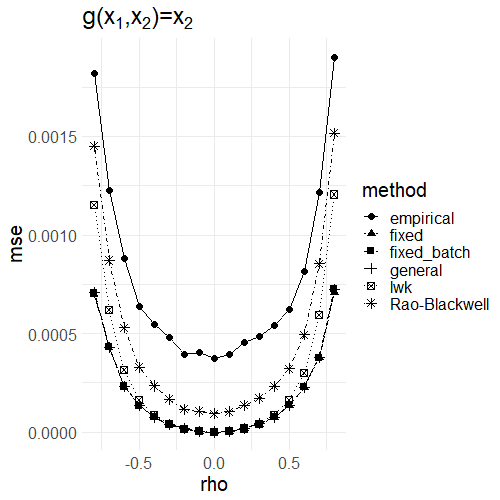}
\caption{}\label{fig:bvn_a}
\end{subfigure}
\begin{subfigure}{0.3\textwidth}
\includegraphics[scale=0.22]{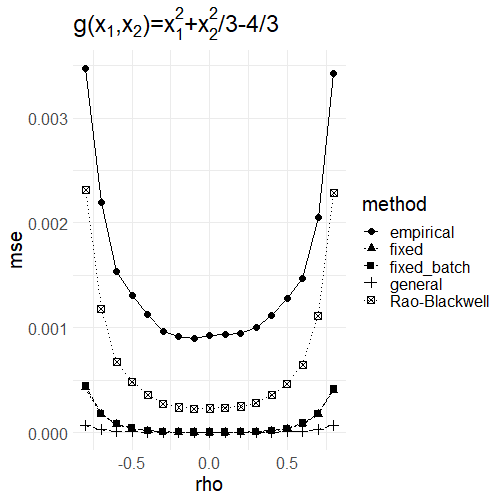}
\caption{}\label{fig:bvn_b}
\end{subfigure}
\begin{subfigure}{0.3\textwidth}
\includegraphics[scale=0.22]{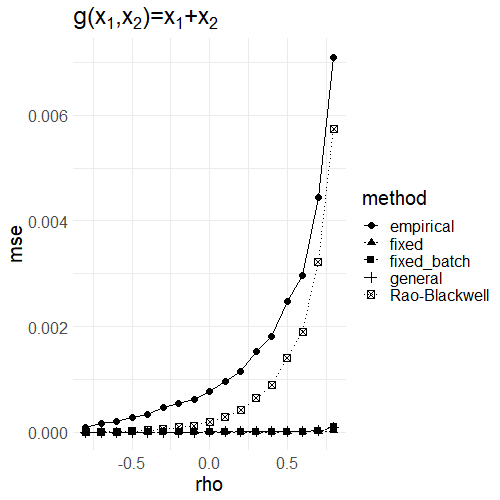}
\caption{}\label{fig:bvn_c}
\end{subfigure}
\caption{Mean squared error (MSE) for different $\rho$ values in a bivariate normal $\pi$ example, with the integrands (\ref{fig:bvn_a}) $g(x_1,x_2)=x_2$, (\ref{fig:bvn_b}) $g(x_1,x_2)=x_1^2+x_2^2/3-4/3$, and (\ref{fig:bvn_c}) $g(x_1,x_2)=x_1+x_2$. The estimator labels general, Rao-Blackwell, empirical, LWK, fixed, and fixed\_batch correspond, respectively, to Equations \eqref{eq:cv1}, \eqref{eq:rb}, \eqref{eq:emp}, \eqref{eq:condEst}, \eqref{eq:cv0} with $\tilde{C}$ estimated by \eqref{eq:cEst} and \eqref{eq:gibbsV}, and \eqref{eq:cv0} with $\tilde{C}$ estimated by \eqref{eq:cEst} and \eqref{eq:genV}.\label{fig:bvn}}
\end{figure}

Our numerical results for the bivariate normal setting are shown in Figure~\ref{fig:bvn}. We compare the simulation mean squared error (MSE) for multiple estimators of $\int \pi(dx)g(x)$, for three different integrands $g$. For each example integrand $g$, we set the control variate basis function $f=g$. For the fixed weight control variate approach, we compare two different estimation methods for the weight $V$ in~\eqref{eq:V}, namely, the estimator~\eqref{eq:gibbsV}, which exploits the simplied representation for $V$ from Corollary~\ref{cor:aVar0}, and the batch means estimator~\eqref{eq:genV}. Empirically, the estimator~\eqref{eq:gibbsV} performs better than the batch means estimator~\eqref{eq:genV}, which is expected since it requires the estimation of fewer covariance terms. We use a batch size of $B=10$. We computed MSE in each setting based on $100$ simulated averages using $M=2000$ draws, at each value of $\rho$.

 Figure~\ref{fig:bvn_a} shows a data augmented setting, where the integrand $g(x_1,x_2)=x_2$, so that $g$ only depends $x_2$. Figure~\ref{fig:bvn_a} compares the simulation asymptotic variances of $S_M/M$ as the bivariate normal correlation coefficient $\rho$ varies for the LWK~\eqref{eq:condEst}, fixed weight control variate~\eqref{eq:cv0}, general control variate~\eqref{eq:cv1}, Rao-Blackwellization~\eqref{eq:rb}, and empirical~\eqref{eq:emp} estimators. We see that the control variate and LWK estimators outperform the empirical and Rao-Blackwell estimators, with the empirical estimator performing the worst. The LWK and control variate estimators perform similarly, although for large $|\rho|$, the control variate estimates outperform the LWK estimates. For $\rho=0$, the LWK estimate is exactly $\Pi_1g(x_1,x_2)=0$ for all $(x_1,x_2)$. Thus, at $\rho=0$, the finite sample performance of the LWK estimate is better than the control variate estimate estimates, which accrue some error in finite samples due the estimation of $\tilde{C}$. This error vanishes asymptotically with $M^{1/2}$ normalization.

Figure~\ref{fig:bvn_b} shows results for the integrand $g(x_1,x_2)=x_1^2+x_2^2/3-4/3$. Since this $g$ depends on both $x_1$ and $x_2$, the approach by~\citet{wongKongLiu} no longer applies. Figure~\ref{fig:bvn_b} compares the variances of the fixed weight and general control variate estimators~\eqref{eq:cv0} and~\eqref{eq:cv1}, as well as the Rao-Blackwellized and empirical estimators~\eqref{eq:rb} and~\eqref{eq:emp}, as $\rho$ varies. For this example, the general control variate estimates outperform the fixed weight estimates. The fixed weight estimates substantially outperform the empirical and Rao-Blackwellized estimates. 

Figure~\ref{fig:bvn_c} shows results for the integrand $g(x_1,x_2)=x_1+x_2$. The control variate estimates (both fixed-weight and general) attain $0$ asymptotic variance, even though the empirical and Rao-Blackwellization estimates have positive asymptotic variance. This can be explained as follows. Taking the random sweep kernel $Q=(\Pi_1+\Pi_2)/2$, we have $Qg
=(1+\rho)g/2$, so that $g$ is an eigenfunction of $Q$ with eigenvalue $\lambda=(1+\rho)/2$. Therefore, taking $c=1/(1-\lambda)$ and $f=g$ gives $g(x)-c\{g(x)-Qg(x)\}=0$ a.e. $\pi$. Thus, the optimal random sweep Gibbs sampling scheme from Theorem~\ref{thm:detVrev} has an asymptotic variance of 0. From Theorem~\ref{thm:detVrev}, we have that the optimal fixed weight deterministic sweep control variate scheme must also attain $0$ asymptotic variance. Figure~\ref{fig:bvn_c} demonstrates that the control variate estimates indeed achieve 0 asymptotic variance, as the MSE for the control variate estimates are nearly exactly 0 except for large $\rho$, where finite sample error in estimating $\tilde{C}$ causes the MSE to be just barely above 0. On the other hand, the empirical and Rao-Blackwell estimators perform much worse, particularly for larger $\rho$.

\subsection{Ising model}

Next, we consider the one-parameter Ising model on a square, $n\times n$ grid of cells. We take $n=20$. The state space is $X=\{-1,1\}^{n^2}$, where $\mathscr{X}=2^X$ is the set of all subsets of $X$. The probability density function of the desired stationary measure with respect to counting measure on $(X,\mathscr{X})$ is $\pi(x)=\exp\{\eta T(x)-\xi(\eta)\}$, where the sufficient statistic $T(x)=\sum_{i\sim j,i<j}x_ix_j$, and the notation $i\sim j$ indicates that $i$ is a neighbor of $j$. Thus, the contribution from a given $i$,$j$ pair with $i\sim j$ is positive when $x_i$ and $x_j$ are equal, and negative otherwise. The term $\xi(\eta)=\log[\sum_{x\in X}\exp\{\eta T(x)\}]$ is a log normalizing constant. We write $x_{-i}$ for the values at all sites of $x$ except site $i$. Also, we use $x^{i}$ to denote the configuration $x$ with the $i$th value flipped, so that $x_i^{i}=-x_{i}$, and $(x^{i})_{-i}=x_{-i}$.

We consider deterministic sweep Gibbs samplers, as well as deterministic sweeps composed of Metropolis-type updates. We first define Gibbs sitewise kernels for each $i=1,...,n$ by \begin{align*}
&\Pi_i(x,\{x'\})=I(x'_{-i}=x_{-i})\pi(x')\{\pi(x)+\pi(x^i)\}^{-1}\\&=I(x'_{-i}=x_{-i})\exp\{\eta T(x')\}\left[\sum_{\substack{x^*\in X:\\x^*_i\in\{-1,1\}\\x^*_{-i}=x_{-i}}}\exp\{\eta T(x^*)\}\right]^{-1}\;\;\;\;\forall x,x'\in X.
\end{align*}Each $\Pi_i$ is a Gibbs kernel with respect to $(h_i,\pi)$ for the coordinate projection $h_i:X\to\mathbb{R},x\to x_{-i}$. Next, we define sitewise Metropolis kernels $Q_i$ by \begin{align*}
	Q_i(x,x')=I(x'=x^{i})0.9a_i(x)+I(x'=x)[0.1+0.9\{1-a_i(x)\}]
	\end{align*} where $a_i(x)=\min\{\pi(x^i)/\pi(x),1\}$. Each Metropolis kernel $Q_i$ corresponds to proposing to flip the value at the $i$th coordinate with probability $0.9$, and then accepting any flip with probability $a_i(x)$. Note that we do not always propose to flip the value at site $i$. It is straightforward to show that the $Q_i(x,x')$ satisfy the reversibility condition~\ref{cond:piReversible}.

	For each of the Gibbs and Metropolis sitewise update types, we consider two different types of compositions of the sitewise updates, so that in total, four Markov chain schemes are considered. First, in the raster sweep, we construct Markov chains $\{X_t\}_{t=0}^{\infty}$ using the update $\Pi_{\sigma^t(1)}$ (resp., $Q_{\sigma^t(1)})$ at each time step $t$, where the sites are traversed sequentially proceeding first down each column of the grid, and then across the columns in order.
	
\begin{figure}[H]
\includegraphics[scale=0.2]{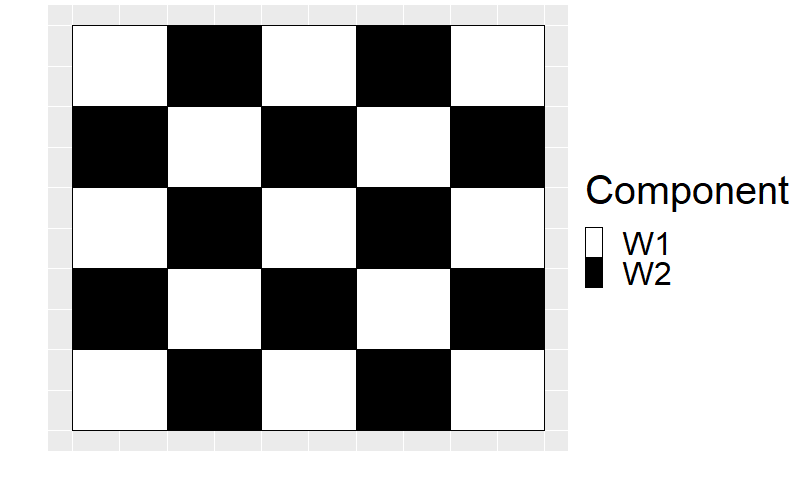}
\caption{Conditional independence structure for the square lattice Ising model. Cells in one component are conditionally independent of each other, given the values at the cells in the other component.}\label{fig:checkerboard}
\end{figure}

	We next consider a checkerboard sweep, where we partition the bipartite lattice into two components $W_1$ and $W_2$, as in Figure~\ref{fig:checkerboard}, and then update each component in sequence. To update each component, we construct composition kernels \begin{align*}
	H_{k}(x,x')=\left\{\prod_{i\notin W_k}\Pi_i\right\}(x,x')\;\;\;\;k=1,2
\end{align*} for the Gibbs kernels and \begin{align*}
	J_{k}(x,x')=\left\{\prod_{i\notin W_k}Q_i\right\}(x,x')\;\;\;\;k=1,2
\end{align*} for the Metropolis kernels. Because of the lattice neighborhood structure of the sites, any ordering of the $\Pi_i$ (resp., $Q_i$) in the compositions $H_k$ (resp., $J_k$) leads to an equivalent transition kernel, and both composition kernels can be implemented using independent Bernoulli draws at every site not in $W_k$. For the checkerboard sweeps, we construct Markov chains $\{X_{t}\}_{t=0}^{\infty}$ by applying the kernel $H_{\sigma(t)}$ (resp., $J_{\sigma(t)}$) at each time $t$. For example, for the Gibbs-based sampler, $H_1$ is used to obtain $X_1$ from $X_0$, and $H_2$ is used to obtain $X_2$ from $X_1$. Thus, at each step, all of the cells are updated in one of the components $W_k$.

It is straightforward to verify that $H_k$ itself is a Gibbs kernel with respect to $(h_{k},\pi)$, where $h_{k}:X\to \mathbb{R}^{|W_k|},x\to x_{W_{k}}$ denotes the coordinate projection which obtains the values in $W_{k}$.

All four sweeps are irreducible with respect to the uniform probability measure on $(X,\mathscr{X})$. Additionally, taking $C=X$, $b=1$, $V_k(x)$ to be the constant function $V_k(x)=2\max_{x'\in X}|T(x')|$ for $k=1,2$, and $g(x)=T(x)-\int\pi(dx')T(x')$ ensures~\ref{cond:drift} holds for chains composed of either the Gibbs and Metropolis updates. Finally, the aperiodicity condition in Assumption~\ref{cond:aperiodicity} holds for the composition chains $P_k^K$ (where $K=n^2$ for the raster sweeps and $K=2$ for the checkerboard sweeps). For both the Gibbs and Metropolis updates, this follows from the fact that \begin{align*}
P_k^K(x,\{x\})>0,\;\;\;\;\forall x\in X
\end{align*} for raster and checkerboard scans with either Gibbs or Metropolis sitewise updates. For the Gibbs sampler chains, we could alternatively have verified~\ref{cond:aperiodicity} by using Lemma~\ref{lem:aperiodic}.

For each sweep type, we considered estimators based on 2000 cycles through the grid. For the raster sweeps, we used a batchsize of $5(400)=2000$, which corresponds to lags resulting from 5 complete cycles through the grid. For the checkerboard sweeps, we used a batchsize of $5(2)=10$, which also corresponds to 5 cycles through the grid.

\begin{figure}[H]
\begin{subfigure}[c]{1\textwidth}
	\centering
\includegraphics[scale=0.2]{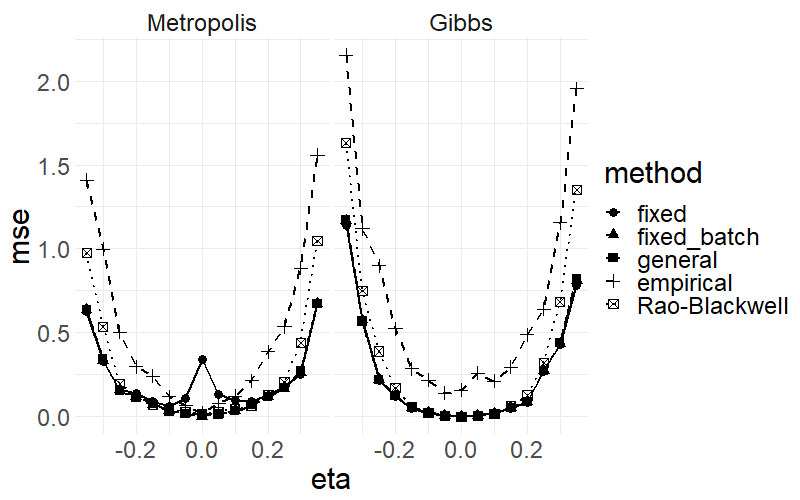}
\caption{Ising model with checkerboard sweep}\label{fig:isingMSE_checkerboard}
\end{subfigure}

\begin{subfigure}[c]{1\textwidth}
	\centering
	\includegraphics[scale=0.2]{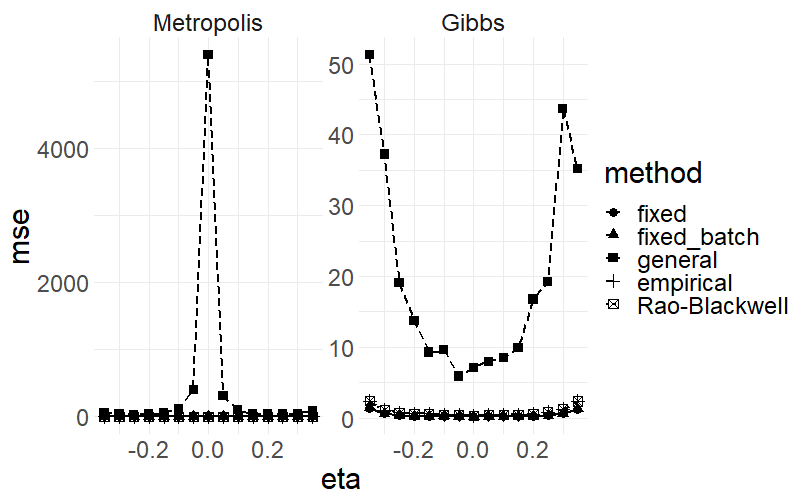}
	\caption{Ising model with raster sweep}\label{fig:isingMSE_raster}
\end{subfigure}

\begin{subfigure}[c]{1\textwidth}
	\centering
	\includegraphics[scale=0.2]{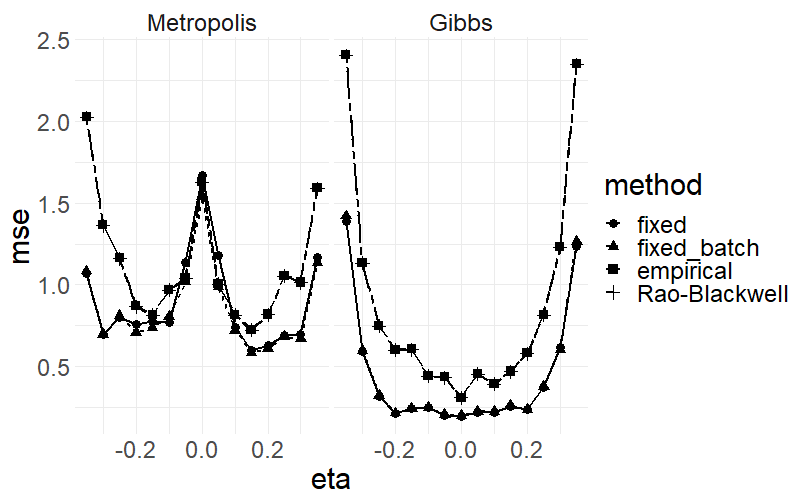}
	\caption{Ising model with raster sweep, with general control variate estimator exluded}\label{fig:isingMSE_rasterII}
\end{subfigure}

\caption{Mean squared error (MSE) for the Ising model simulation example at different values of $\eta$, for deterministic raster and checkerboard sweeps.}\label{fig:isingMSE}
\end{figure}

\begin{figure}[H]
	\begin{subfigure}[c]{1\textwidth}
	\centering
	\includegraphics[scale=0.2]{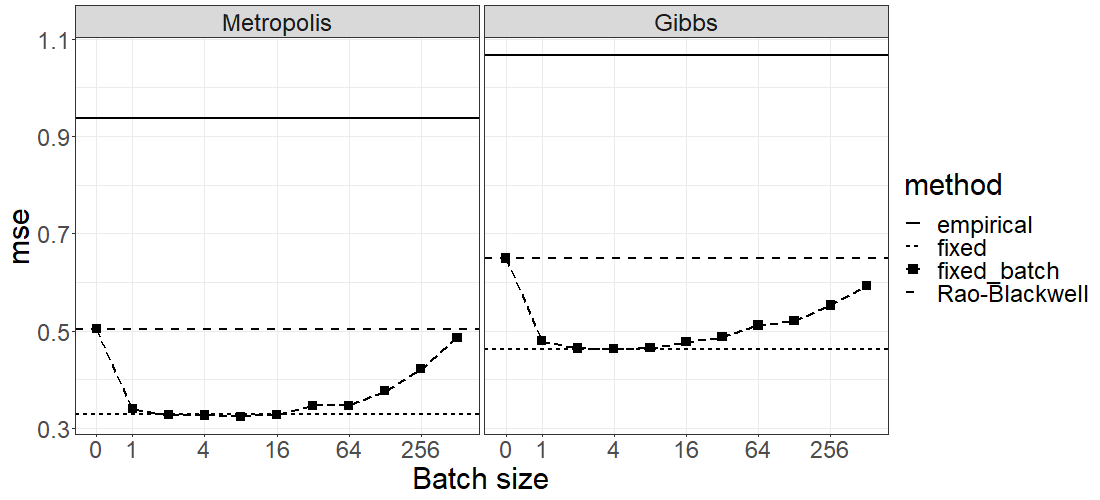}
	\caption{Checkerboard sweep}\label{fig:isingBatch_checkerboard}
	\end{subfigure}
	
	\begin{subfigure}[c]{1\textwidth}
		\centering
		\includegraphics[scale=0.2]{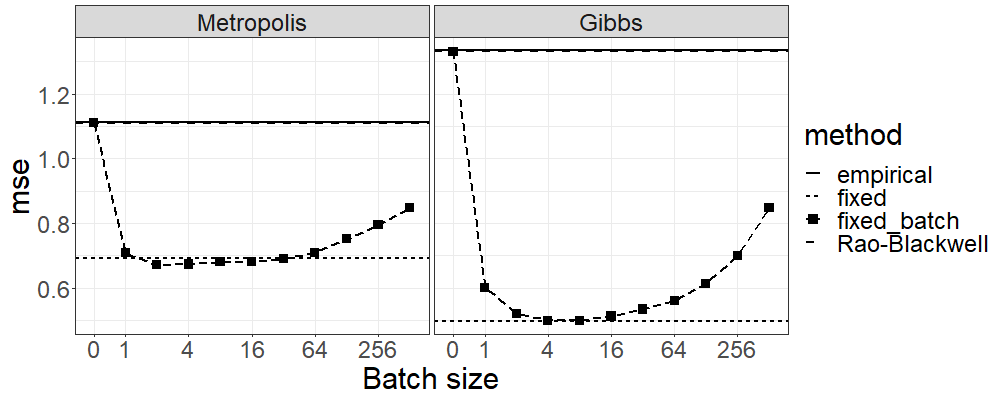}
		\caption{Raster sweep}\label{fig:isingBatch_raster}
	\end{subfigure}

	\caption{Mean squared error (MSE) for the fixed weight control variate method with checkerboard and raster sweeps, for the Ising model with $\eta=0.3$, based on 100 simulated means at each batch size $B$, where each simulated mean used $M=2000$ draws.}	\label{fig:isingBatch}
\end{figure}

In Figures~\ref{fig:isingMSE} and~\ref{fig:isingBatch}, method=``fixed'' and ``fixed\textunderscore batch'' correspond to estimating $\tilde{C}$ via~\eqref{eq:cEst} with $V$ estimated via \eqref{eq:gibbsV} and~\eqref{eq:genV}, respectively. In Figure~\ref{fig:isingMSE}, mean squared error (MSE) is computed based on the empirical average squared error of 100 estimated means, where each estimated mean used 2000 complete cycles through the grid. That is, $M=2000(n^2)=800000$ for the raster sweeps, and $M=2000(2)=4000$ for the checkerboard sweeps.

Figure~\ref{fig:isingMSE} shows the MSE for the checkerboard and raster sweeps with Metropolis and Gibbs updates. Figure~\ref{fig:isingMSE_checkerboard} shows the performance of the estimators for checkerboard sweep. The control variate estimates all perform well for the Gibbs sampler. The Metropolis based control variate estimates perform well except for the fixed weight control variate estimate. This is expected since the fixed weight estimate without batch means uses the formula~\eqref{eq:gibbsV}, which is exact for Gibbs sampling but not for Metropolis sampling. The Metropolis control variate schemes perform well when using the proper batch means formulas to estimate the control variate weight.

Figure~\ref{fig:isingMSE_raster} shows that for the raster sweep, the general control variate estimator performs much worse than the other estimators, likely due to the fact that $n^2$ weights must be estimated for this scheme. 

Figure~\ref{fig:isingMSE_rasterII} shows the MSE for the raster sweep, for the estimators remaining after excluding the general control variate estimate. In Figure~\ref{fig:isingMSE_rasterII}, the empirical and Rao-Blackwellized schemes nearly overlap for both Gibbs and Metropolis schemes. For Gibbs sampling, the fixed weight estimator based on~\eqref{eq:gibbsV} performs best, as expected, but the fixed weight batch estimator also performs well. For Metropolis sampling, the fixed weight estimators perform similarly, but the MSE for the batch means estimator is often smaller than for the fixed weight estimator based on~\eqref{eq:gibbsV}.

For each value of $\eta$, we estimated the true value of $\int\pi(dx)g(x)$ using a long checkerboard sweep run with 100000 complete cycles through the grid, so that $M=100000(2)=200000$. We used the Rao-Blackwellized estimator to compute the means.

We also examined the effect of differing batch sizes on the performance of the various estimates (Figure~\ref{fig:isingBatch}). Our study allowed us to confirm three notable theoretical predictions about the performance of the control variate estimators. First, for Gibbs samplers, the optimal fixed weight control variate formula based on~\eqref{eq:gibbsV} (horizontal line with smaller dashes) always performed better than the corresponding batch means approach, as expected. For the Gibbs samplers, the batch means estimators performed best near a batch size of $B=4$. At this batch size, the estimation performances were nearly identical to, but slightly worse than, the estimates using~\eqref{eq:gibbsV}. Second, our results for each setting demonstrate empirically that setting the batch size $B=0$ is asymptotically equivalent to the Rao-Blackwellization approach~\eqref{eq:rb} (horizontal line with larger dashes). Third, for the Metropolis samplers, using the optimal batch size in each setting leads to a better control variate weight than using the fixed weight estimator~\eqref{eq:gibbsV} (horizontal line with smaller dashes).

\section{Conclusions and discussion\label{sec:conclusion}}

The results here can be used to improve the efficiency of Bayesian statistical analyses that involve deterministic sweep Gibbs samplers, or hybrid Gibbs-Metropolis sweep samplers. This work also provides a rigorous foundation for the use of Rao-Blackwellization in the deterministic sweep Gibbs sampling context, where previously the main justifications for the procedure have come through heuristic uses of Jensen's inequality or the classical Rao-Blackwell Theorem for independent data. 

Future work may include investigating good choices of control variate basis functions. In addition, we believe our results will be directly relevant in stochastic approximation settings. Typically in these settings, an iterate sequence $\theta_t$ will converge to some optimum $\theta^*$ at a rate such that a distributional convergence result \begin{align*}
t^{1/2}(\theta_t-\theta^*)\overset{d}{\to}N\{0,J(\theta^*)^{-1}V(\theta^*)J(\theta^*)^{-1}\}
\end{align*} holds~\citep[see, e.g.,][]{Benveniste,kushnerYin,fort2015central}. The asymptotic variance of the iterates $J(\theta^*)^{-1}V(\theta^*)J(\theta^*)^{-1}$ essentially sandwiches the asymptotic variance of the gradient estimates $V(\theta^*)$ near $\theta^*$ between inverse Hessian terms $J(\theta^*)^{-1}$. Thus, using our approaches to reduce $V(\theta^*)$ can plausibly lead to reductions in the iterate variance. Variance reduction methods are already used in many implementations of stochastic approximation algorithms, although mostly in the Rao-Blackwellization sense of Corollary~\ref{cor:naiveVec} rather than the control variate approaches discussed here, which may be superior. For example, in the machine learning community, deep Boltzmann machines~\citep{hinton2010} have been used to model datasets containing complex structure such as the MNIST handwritten digit dataset, and Rao-Blackwellization for variance reduction has been used to make the model fitting process more efficient for these algorithms. 

    \putbib
\end{bibunit}
\newpage
\begin{bibunit}

\begin{center}
    {\Large\bf SUPPLEMENTARY MATERIAL}
\end{center}

\appendix
The proofs of Lemmas~\ref{lem:periodic}--\ref{lem:exampleDrift} are in Appendix~\ref{app:appA}. Appendix~\ref{app:appB} contains (in order), the proofs of Theorem~\ref{thm:vecG}, Corollaries~\ref{cor:aVar0}--\ref{cor:naiveVec}, Theorem~\ref{thm:detVrev}, and Propositions~\ref{prop:poisEq}--\ref{prop:lwkProp}.

\section{Proofs of Lemmas\label{app:appA}}

In Lemmas~\ref{lem:periodic}--\ref{lem:aperiodic}, we show that the length $K$ composition kernels $P_k^{K}$ are aperiodic when the $\Pi_k$ are Gibbs kernels for $k=1,...,K$.

\begin{lem}\label{lem:periodic} Assume~\ref{cond:projection}. Suppose $P_k^{t}I_A(x)=I_A(x)$ a.e. $\pi$ for some $t>1$. Then $P_{\sigma(k)}^{t-1}I_A=I_A$ a.e. $\pi$.
\end{lem}

Proof:
We have \begin{align*}&\braket{I_A,I_A}=\braket{P_k^{t}I_A,P_k^{t}I_A}=\braket{P_{\sigma(k)}^{t-1}I_A,P_k^{t}I_A}\\&\leq \braket{P_{\sigma(k)}^{t-1}I_A,P_{\sigma(k)}^{t-1}I_A}^{1/2}\braket{P_k^{t}I_A,P_k^{t}I_A}^{1/2}\\&=\braket{P_{\sigma(k)}^{t-1}I_A,P_{\sigma(k)}^{t-1}I_A}^{1/2}\braket{I_A,I_A}^{1/2}\end{align*} where the first and final equalities follow since $P_k^{t}I_A=I_A$ a.e. $\pi$, the second equality follows from reversibility and idempotence of $\Pi_k$, and the inequality follows from the Cauchy-Schwarz inequality. Jensen's inequality gives $\braket{P_{\sigma(k)}^{t-1}I_A,P_{\sigma(k)}^{t-1}I_A}\leq \braket{I_A,I_A}$, which from the preceding implies $\braket{P_{\sigma(k)}^{t-1}I_A,P_{\sigma(k)}^{t-1}I_A}=\braket{I_A,I_A}$. Since $\braket{P_{\sigma(k)}^{t-1}I_A,P_{\sigma(k)}^{t-1}I_A}=\braket{I_A,I_A}$, applying the Cauchy-Schwarz inequality to $P_{\sigma(k)}^{t-1}I_A$ and $P_k^{t}I_A$ implies $P_{\sigma(k)}^{t-1}I_A(x)=P_k^{t}I_A(x)=I_A(x)$ a.e. $\pi$.\qed

Lemma~\ref{lem:pitoPsi} below relates the stationary measure $\pi$ to the irreducibility measure $\psi$.
\begin{lem}\label{lem:pitoPsi}
	Assume~\ref{cond:piStationary} and~\ref{cond:irreducibility}, and suppose $\psi(A)>0$ for some $A\in\Xs$, and $P_k^K(x,A)=1$ for all $x\in A$. Then $\pi(A)=1$.
\end{lem}

Proof of Lemma~\ref{lem:pitoPsi}: Note that for any $A\in\Xs$, \begin{align*}
	&\pi(A)=\sum_{t=1}^{\infty}2^{-t}\pi(A)=\sum_{t=1}^{\infty}2^{-t}\int\pi(dx)P_k^{Kt}I_A(x)\nonumber
\end{align*}

Now, suppose some set $A\in\Xs$ satisfies $\psi(A)>0$ and $P_k^K(x,A)=1$ for all $x\in A$. Then \begin{align*}
	&\pi(A)=\sum_{t=1}^{\infty}2^{-t}\int\pi(dx)P_k^{Kt}I_A(x)\nonumber\\
	&=\sum_{t=1}^{\infty}2^{-t}\int\pi(dx)I_AP_k^{Kt}I_A(x)\\
	&+\sum_{t=1}^{\infty}2^{-t}\int\pi(dx)I_{A^C}P_k^{Kt}I_A(x)\\
	&=\pi(A)+\int\pi(dx)I_{A^C}\sum_{t=1}^{\infty}2^{-t}P_k^{Kt}I_A(x)
\end{align*} Thus, $I_{A^C}\sum_{t=1}^{\infty}2^{-t}P_k^{Kt}I_A(x)=0$ a.e. $\pi$. But from $\psi$-irreducibility of $P_k^K$, the infinite sum is positive for all $x$. This implies $I_{A^C}=0$ a.e. $\pi$, and thus $\pi(A)=1$.\qed

Now, we finish the proof of aperiodicity.
\begin{lem}\label{lem:aperiodic} Under~\ref{cond:projection} and~\ref{cond:irreducibility}, the transition kernels $P_k^K$ are aperiodic.
\end{lem} Proof: Consider some arbitrary $k$ from $1,...,K$. From Theorem 5.4.4 of \citet{meyn2009markov} and the $\psi$-irreducibility of $P_k^{K}$, there exists an integer $d$ and a collection of sets $D_1,...,D_d$ satisfying \begin{enumerate}
	\item $P_k^K(x,D_{i+1})=1$ for $x\in D_i$, $i\equiv 0,...,d-1$ (mod $d$)
	\item $\psi\{(\cup_{i=1}^{d}D_i)^C\}=0$
	\item $D_1,...,D_{d}$ are disjoint
\end{enumerate}

We show that $P_k^K$ is aperiodic by showing that $d=1$ is the largest integer such that 1-3 hold for a collection of sets $D_1,...,D_d$. Suppose to the contrary that $d>1$ for a collection of sets $D_1,...,D_d$ satisfying 1-3. From Lemma~\ref{lem:pitoPsi}, we have $\pi(\cup_{i=1}^{d}D_i)=1$. Thus, $P_k^{Kd}I_{D_i}=I_{D_i}$ a.e. $\pi$ for each $i$. Now, $K(d-1)$ applications of Lemma~\ref{lem:periodic} imply $P_k^{K}I_{D_i}=I_{D_i}$ a.e. $\pi$ for each $i$. Additionally, $\pi(D_i)>0$ for at least one $i$, so for this $i$, $H=\{x\in D_i:P_k^KI_{D_i}=1\}$ is non-empty. But this is a contradiction, since $P_k^K(x,D_{i+1})=1$ for all $x\in H$. Thus, any collection of sets $D_1,...,D_d$ satisfying 1-3 must have $d=1$. This proves the result.\qed



\begin{lem}\label{lem:asConst}
	Assume~\ref{cond:piStationary}--\ref{cond:aperiodicity} hold. Suppose $\Pi_kf=f$ a.e. $\pi$ for each $k=1,...,K$, where $f:(X,\mathscr{X})\to (\mathbb{R},\mathscr{R})$. Then $f$ is constant a.e. $\pi$.
\end{lem}

Proof: Suppose $\Pi_kI_A=I_A$ a.e. $\pi$ for each $k$, for some $A\in\mathscr{X}$. We first show $\pi(A)=0$ or 1. Define a Markov chain $\{Y_t\}_{t=0}^{\infty}$ by the initial distribution $Y_0\sim \pi$ and the transition kernel $P_1^{K}$. Then we have $I_A(Y_t)=I_A(Y_0)$ almost surely for all $t$. Now under~\ref{cond:piStationary}--\ref{cond:aperiodicity}, a  Law of Large Numbers holds, so that $M^{-1}\sum_{t=0}^{M-1}I_A(Y_t)\overset{a.s.}{\to}\int \pi(dx)I_A(x)$. Thus, $I_A(Y_0)=\int\pi(dx)I_A(x)$ almost surely. This implies $\int \pi(dx)I_A(x)=0$ or $1$.

Now, let $\mathscr{H}=\{A\in\mathscr{F}:\Pi_kI_A=I_A\;a.e.\;\pi,\;\forall k\}$. We show $\mathscr{H}$ is a $\sigma$-field. The empty set $\phi\in \mathscr{H}$ and the state space $X\in\mathscr{H}$. Also, for any arbitrary $A\in\mathscr{H}$, we have $\pi(A)=0$ or $\pi(A)=1$. This implies $A^C\in \mathscr{H}$ for each $A\in\mathscr{H}$. Finally, for $\{A_n\}_{n=1}^{\infty}$ with each $A_n\in\mathscr{H}$, we have $\pi(\cup_{n=1}^{\infty}A_n)=0$ or 1, so that $\cup_{n=1}^{\infty}A_n\in\mathscr{H}$. Thus, $\mathscr{H}$ is a $\sigma$-field.

We have $\Pi_kf=f$ a.e. $\pi$ by assumption. We now show $f$ is $\mathscr{H}$-measurable. Let $B\in\mathscr{R}$ given. Define $A=f^{-1}(B)$. Also, define $A_k=(\Pi_kf)^{-1}(B)$ for each $k$. We have $I_A=I_{A_k}$ a.e. $\pi$ for each $k$ since $f=\Pi_kf$ a.e. $\pi$ for each $k$. Thus $A=f^{-1}(B)\in \mathscr{H}$. Since $B$ was arbitrary, $f$ is $\mathscr{H}$ measurable.

Finally, we show that $f$ is constant a.e. $\pi$. Without loss of generality, we assume $f\geq 0$. For general $f$, we may use the standard decomposition of $f$ into positive and negative components and apply the following reasoning to each component \citep[see, e.g., Chapter 1 of][]{shao2003mathematical}. Since $f$ is $\mathscr{H}$-measurable, we can construct a sequence $\{f_n\}_{n=1}^{\infty}$ of $\mathscr{H}$-measurable simple functions such that $f_n\uparrow f$ pointwise. Now, let $b=\int \pi(dx)f(x)$. Then $\pi(\{x:f(x)>b\})=\underset{n\to\infty}{\lim}\pi(\{x:f_n>b\})=0$, where we used $\{f_{n}>b\}\subset\{f_{n+1}>b\}$ for each $n$ by monotone convergence of the $f_n$, as well as the fact that $f_n$ are constant a.e. $\pi$ so that $\int \pi(dx)f_{n}\leq b$ implies $f_n\leq b$ a.e. $\pi$. Thus, $f(x)=b$ a.e. $\pi$. This completes the proof.\qed

\begin{lem}\label{lem:smallCondExp}
Assume~\ref{cond:piStationary}--\ref{cond:aperiodicity} hold, and suppose~\ref{cond:smallf} holds for the function $f:X\to\mathbb{R}^p$. Then for any $a\in\mathbb{R}^p$, we have \begin{align*}a^T\left[\sum_{k}\int\pi(dx)\{ff^T-(\Pi_kf)(\Pi_kf)^T\}\right]a=0\end{align*} if and only if $a^Tf=b$ a.e. $\pi$ for some constant $b$.
\end{lem}

Proof: Follows easily from Lemma~\ref{lem:asConst}.\qed

\begin{lem}\label{lem:completeSquare}
	Let $U$ be a symmetric positive semidefinite $p\times p$ matrix, and $V$ be a $p\times d$ matrix such that $a^TV=0_{1\times d}$ for any vector $a$ with $Ua=0$. Then for all $p$-vectors $a$, \begin{align*}
		U^{\dagger}V\in\underset{C\in\mathbb{R}^{p\times d}}{\arg\min} \;a^T\{C^TUC-C^TV-V^TC\}a
	\end{align*}
\end{lem}

Proof: First, we note that when $U=0_{p\times p}$, then $V=0$ also, so that $a^T\{C^TUC-C^TV-V^TC\}a=0$ for all $a$, for any choice of $C$. In particular, $a^T\{C^TUC-C^TV-V^TC\}a=0$ for all $a$ for $C=U^{\dagger}V$. Otherwise, since $U$ is symmetric positive semidefinite, we may write $U=QDQ^T$ where $Q$ is a $p\times r$ matrix with $r\leq p$ orthonormal columns, and $D$ is a $r\times r$ diagonal matrix with strictly positive diagonal entries. Further, for any $C\in\mathbb{R}^{p\times d}$, we may write $C=QR+B$ where $R\in\mathbb{R}^{r\times d}$, $B\in\mathbb{R}^{p\times d}$, and $Q^TB=0_{r\times d}$.

It can be checked that the value of the $B$ component of $C$ does not affect the value of $a^T(C^TUC-C^TV-V^TC)a$, so that minimizers of the form $C=QR$ exist. When $C=QR$, we have
\begin{align}\label{eq:completeSquare}
	C^TUC-C^TV-V^TC=X^TX-V^TQD^{-1}Q^TV
\end{align} where $X=D^{1/2}Q^TC-D^{-1/2}Q^TV$. The second term in \eqref{eq:completeSquare} does not depend on $C$. Now, $a^TX^TXa\geq 0$ for arbitrary $C$. But $U^{\dagger}=QD^{-1}Q^T$, so that taking $C=U^{\dagger}V$ gives $X=D^{1/2}Q^TQD^{-1}Q^TV-D^{-1/2}Q^TV=0$. Thus, $a^T\{C^TUC-C^TV-V^TC\}a$ is minimized for each $a$ whenever $C=U^{\dagger}V$. This completes the proof. \qed





In Lemma~\ref{lem:Qlem}--\ref{lem:Qlem1}, we take $K=2$ and $Q=(\Pi_1+\Pi_2)/2$.
\begin{lem}\label{lem:Qlem}
	Assume~\ref{cond:projection}--\ref{cond:drift} and~\ref{cond:integrability}--\ref{cond:smallg}. Then we have
		$\sum_{t=1}^{\infty}|Q^tg(x)|$ is square integrable with respect to $\pi$, and $\sum_{t=1}^{\infty}Q^tg(x)=\sum_{t=1}^{\infty}(P_1^t+P_2^t)g(x)$ a.e. $\pi$.
\end{lem}
Proof: We show the result for scalar $g:X\to\mathbb{R}$. The result follows for general $g:X\to\mathbb{R}^p$ by applying the reasoning below elementwise.

Note $Q^tg=\{(\Pi_1+\Pi_2)/2\}^tg=\sum_{i=1}^{t}2^{-t}\binom{t-1}{i-1}(P_1^i+P_2^i)g$ a.e. $\pi$ via idempotence of the $\Pi_k$, so that \begin{align}
	&\sum_{n=1}^{M}|Q^ng(x)|\leq \sum_{n=1}^{M}2^{-n}\sum_{i=1}^{n}\binom{n-1}{i-1}\{|P_1^ig(x)|+|P_2^ig(x)|\}\nonumber\\
	&=\sum_{i=1}^{M}\sum_{n=1}^{M}\{|P_1^ig(x)|+|P_2^ig(x)|\}2^{-n}\binom{n-1}{i-1}I(i\leq n)\label{eq:binomial}\\
	&\leq \sum_{n=1}^{M}(|P_1^ng(x)|+|P_2^ng(x)|)\nonumber
\end{align} a.e. $\pi$, where the second inequality follows because \begin{align*}
	\sum_{n=i}^{\infty}2^{-n}\binom{n-1}{i-1}=\sum_{r=1}^{\infty}2^{-(i-1+r)}\binom{i-2+r}{i-1}=1,
\end{align*} which itself is a well-known identity related to the pdf of a negative binomial random variable.

The Assumptions~\ref{cond:projection}--\ref{cond:drift} imply~\ref{cond:piStationary}--\ref{cond:aperiodicity} hold, and since~\ref{cond:integrability}--\ref{cond:smallg} also hold, we have $\sum_{t=1}^{\infty}|P_1^tg(x)|+|P_2^tg(x)|<\infty$ a.e. $\pi$ by Proposition~\ref{prop:poisEq}. Thus $\sum_{n=1}^{\infty}|Q^ng(x)|$ converges a.e. $\pi$, and \begin{align*}
	&\sum_{n=1}^{\infty}Q^ng(x)=\lim_{M\to\infty}\sum_{i=1}^{M}\sum_{n=1}^{M}\{P_1^ig(x)+P_2^ig(x)\}2^{-i}\binom{n-1}{i-1}I(i\leq n)\\
	&=\sum_{n=1}^{\infty}P_1^ng(x)+P_2^ng(x),
\end{align*} a.e. $\pi$, where the final equality follows from using Fubini's Theorem and the identity $\sum_{i=n}^{\infty}2^{-i}\binom{i-1}{n-1}=1$.\qed

\begin{lem}\label{lem:Qlem1}
	Assume~\ref{cond:projection},~\ref{cond:irreducibility}--\ref{cond:drift}, and~\ref{cond:integrability}--\ref{cond:smallf} hold. Assume $\{X_t\}_{t=0}^{\infty}$ is defined as in Theorem~\ref{thm:vecG}. Take $S_M$ as in~\eqref{eq:cv0} and define $H_M=\sum_{t=0}^{M-1}g(X_t)-C^T\{f(X_t)-Qf(X_t)\}$. Then $M^{-1/2}(S_M-H_M)\overset{a.s.}{\to}0$.
\end{lem}

Proof of Lemma~\ref{lem:Qlem1}: We prove the result assuming $X_0\sim \pi$. When $X_0\sim \nu$ for general $\nu$, the result can be shown with a coupling argument as sketched in Theorem~\ref{thm:vecG}. 

We prove the result for scalar $f,g:X\to\mathbb{R}$ and $C=c\in\mathbb{R}$. The extension to $g:X\to\mathbb{R}^d$, $f:X\to\mathbb{R}^p$, and $C\in\mathbb{R}^{p\times d}$ can be shown by applying the univariate result elementwise for each of the $d$ elements of $M^{-1/2}(S_M-H_M)$.

For univariate $f,g$ and scalar $c$, we have \begin{align*}S_M-H_M&=\sum_{t=0}^{\floor{(M-1)/2}}c\Pi_1f(X_{2t})+\sum_{t=0}^{\floor{(M-2)/2}}c\Pi_2f(X_{2t+1})\\
&-\sum_{t=0}^{M-1}c\{\Pi_1f(X_t)+\Pi_2f(X_t)\}/2.\end{align*} Now, when $t$ is even, $\Pi_1f(X_t)=\Pi_1f(X_{t+1})$ almost surely, and similarly, $\Pi_2f(X_t)=\Pi_2f(X_{t+1})$ almost surely when $t$ is odd. Thus \begin{align*}
	&\sum_{t=0}^{\floor{(M-1)/2}}c\Pi_1f(X_{2t})-\sum_{t=0}^{M-1}c\Pi_1f(X_{t})/2\\
	&=\begin{cases}
		0 & M\text{ even}\\
		c\Pi_1f(X_{M-1})/2 & M\text{ odd}
	\end{cases}
	\end{align*}
almost surely, and
	\begin{align*}
	&\sum_{t=0}^{\floor{(M-2)/2}}c\Pi_2f(X_{2t+1})-\sum_{t=0}^{M-1}c\Pi_2f(X_{t})/2\\
	&=\begin{cases}
		-c\Pi_2f(X_0)/2 & M\text{ odd}\\
		-c\Pi_2f(X_0)/2 + c\Pi_2f(X_{M-1})/2& M\text{ even}\\
	\end{cases}
\end{align*}
	almost surely. Thus, \begin{align*}
		&M^{-1/2}|S_M-H_M|\\
		&\leq M^{-1/2}c\{|\Pi_2f(X_0)/2|+|\Pi_1f(X_{M-1})/2|+|\Pi_2f(X_{M-1})/2|\}\overset{a.s.}{\to}{0}
\end{align*} as $M\to\infty$ by applying the Strong Law of Large Numbers along the $K=2$ subchains to the function $\{(\Pi_1+\Pi_2)f(x)/2\}^2$. \qed

\begin{lem}\label{lem:exampleDrift}
    Consider the Gibbs kernels $\Pi_1$ and $\Pi_2$ from the bivariate normal Gibbs sampling example. Then there exists an $r>0$ such that Assumption~\ref{cond:drift} holds for the composition kernels $P_1^2$ and $P_2^2$ with the choices $V_1(x)=x_1^2+rx_2^2+1$ and $V_2(x)=rx_1^2+x_2^2+1$.
    \end{lem}

    Proof: First, we show that the $\{X_{2n}\}_{n=0}^{\infty}$ and $\{X_{2n+1}\}_{n=1}^\infty$ chains are T-chains, in the sense of \citet{meyn2009markov}. To do this, we show that the composition kernels $P_{1}^2(x,\cdot)$ and $P_2^2(x,\cdot)$ are strong Feller chains. Since the kernels $P_{k}^2(x,\cdot)$ are aperiodic (Lemma \ref{lem:aperiodic}) and $\psi$-irreducible, this will imply from part $(ii)$ of Theorem 6.2.5 in \citet{meyn2009markov} that every compact subset of $\mathbb{R}^2$ is small.

    To show that $P_1^2(x,\cdot)$ is strong Feller, we check that \begin{align*}
        \underset{n}{\lim\inf}P_1^2(x_n,A)\geq P_1^2(x,A)
    \end{align*} for any $A\in\mathscr{R}^2
    $ and sequence $\{x_n\}_{n=1}^{\infty}$ with $x_n=(x_{1n},x_{2n})\in\mathbb{R}^2$ and $x_n\to x$.

    Let $A\in\mathscr{R}^2$, and suppose $\{x_n\}_{n=1}^{\infty}$ is a sequence in $\mathbb{R}^2$ with $x_n\to x^*$. We have \begin{align*}
        &\int P_1^2(x_n,dx')I_A(x')\\
        &=B\int\exp\{-(1-\rho^2)^{-1}(x_2'-\rho x_{1n})^2/2\}\exp\{-(1-\rho^2)^{-1}(x_1'-\rho x_2')^2/2\}I_A(x')dx_1'dx_2'
    \end{align*} where $x'=(x_1',x_2')$, and the constant $B=\{2\pi(1-\rho^2)\}^{-1}$ does not depend on $x_n$ or $x'$. Now, since $x_n\to x^*$, we have in particular that $x_{1n}\to x_1^*$. Thus, by continuity, \begin{align*}
        &\underset{n}{\lim\inf}\exp\{-(1-\rho^2)^{-1}(x_2'-\rho x_{1n})^2/2\}\exp\{-(1-\rho^2)^{-1}(x_1'-\rho x_2')^2/2\}\\
        &=\exp\{-(1-\rho^2)^{-1}(x_2'-\rho x_1^*)^2/2\}\exp\{-(1-\rho^2)^{-1}(x_1'-\rho x_2')^2/2\}
    \end{align*}

    Therefore, from Fatou's Lemma, we have $\underset{n}{\lim\inf}P_1^2(x_n,A)\geq P_1^2(x^*,A)$, so $P_1^{2}$ is strong Feller. The proof that $P_2^2$ is also strong Feller is similar. Thus, all compact sets are small for $P_1^2$ and $P_2^2$ from Theorem 6.2.5 in \citet{meyn2009markov}.

    Now, consider $V_1(x)=x_1^2+rx_2^2+1$ where $0<r<(1-\rho^4)$. Take $\lambda_1=\rho^4+r$. Then $\Pi_2V_1(x)=\rho^2x_2^2+(1-\rho^2)+rx_2^2+1$ and \begin{align*}		\Pi_1\Pi_2V_1(x)&=(\rho^4+r\rho^2)x_1^2+(1+\rho^2+r)(1-\rho^2)+1\\
        &\leq \lambda_1 V_1(x)+(1+\rho^2+r)(1-\rho^2)+1-\lambda_1-\lambda_1rx_2^2
    \end{align*}

    Now, take $b=(1+\rho^2+r)(1-\rho^2)+1-\lambda_1$ and $c>0$ such that $\lambda_1rc^2\geq (1+\rho^2+r)(1-\rho^2)+1-\lambda_1$. Then we have $P_1^2V_1(x)=\Pi_1\Pi_2V_1(x)\leq \lambda_1V_1(x)+bI_C(x)$ where $C=[-c,c]\times [-c,c]$. Since $C$ is compact, $C$ is small, so~\ref{cond:drift} is satisfied for $P_1^2$.

	Similarly, it can be shown that $P_2^2$ also satisfies~\ref{cond:drift}. \qed

\section{Proofs of Theorems, Corollaries, and Propositions\label{app:appB}}

\textbf{Proof of Theorem~\ref{thm:vecG}}:

First, we consider the case when the initial measure $\nu=\pi$. In this case, the law of $\{X_t\}_{t=0}^{\infty}$ is $P_\pi$ and $X_0\sim \pi$. To simplify notation, we will prove the result in the univariate case where $g:X\to\mathbb{R}$, $f:X\to\mathbb{R}$, and $C_k=c_k\in\mathbb{R}$. In the remainder, for notational clarity, we will use the conventions $\Pi_t=\Pi_{\sigma^t(1)}$, $c_t=c_{\sigma^t(1)}$, and $\hat{g}_t=\hat{g}_{\sigma^t(1)}$,
so that
\begin{align}
&S_M=\sum_{t=0}^{M-1}g(X_t)-c_{t}f(X_t)+c_{t+1}\Pi_{t}f(X_t)\label{eq:sn}\\
&=\sum_{t=0}^{M-1}\hat{g}_{t}(X_t)-\Pi_t\hat{g}_{t+1}(X_t)-c_{t}f(X_t)+c_{t+1}\Pi_{t}f(X_t)\nonumber\\
&=\hat{g}_0(X_0)-\hat{g}_M(X_M)-c_0f(X_0)+c_Mf(X_M)\nonumber\\
&+\sum_{t=0}^{M-1}\hat{g}_{t+1}(X_{t+1})-\Pi_{t}\hat{g}_{t+1}(X_t)-c_{t+1}\{f(X_{t+1})-\Pi_{t}f(X_t)\}\nonumber
\end{align} a.e. $P_{\pi}$, where we are using the identity $\hat{g}_t-\Pi_t\hat{g}_{t+1}=g$ a.e. $\pi$ from Proposition~\ref{prop:poisEq} in the second equality, and rearranging the sum in the third equality.

The term $U_{M}:=\sum_{t=0}^{M-1}\hat{g}_{t+1}(X_{t+1})-\Pi_{t}\hat{g}_{t+1}(X_t)-c_{t+1}\{f(X_{t+1})-\Pi_{t}f(X_t)\}$ is an $L_2$ martingale (since $X_0\sim \pi$, and the $\hat{g}_{t}$ are square integrable with respect to $\pi$ from Proposition~\ref{prop:poisEq}). The remainder term $\hat{g}_0(X_0)-\hat{g}_M(X_M)-c_0f(X_0)+c_Mf(X_M)$ will be shown to be small using the Law of Large Numbers for Markov chains. Thus, we expect the asymptotic behavior of $S_M$ to be similar to that of $U_M$, and we will apply a central limit theorem for martingales to deal with this term.

We now introduce a martingale central limit theorem, Theorem~\ref{thm:mtgCLT}, which follows immediately from Theorem 3.2, Corollary 3.1 of \citet{hall1980martingale}. We use $\overset{p}{\to}$ to denote convergence in probability.

\begin{thma}\label{thm:mtgCLT}
Let $\{S_{ni},\mathscr{F}_{ni},1\leq i\leq k_n,n\geq 1\}$ be a zero-mean, square integrable martingale array with differences $Y_{ni}=S_{ni}-S_{n,i-1}$ ($S_{n0}:=0$). Suppose \begin{enumerate}
	\item (conditional Lindeberg) \label{cond:condLind}for all $\epsilon>0$, $\sum_{i=1}^{k_n}E\{Y_{ni}^2I(|Y_{ni}|>\epsilon)|\mathscr{F}_{n,i-1}\}\overset{p}{\to}0$
	\item (converging conditional variances) \label{cond:convVar}$\sum_{i=1}^{k_n}E(Y_{ni}^2|\mathscr{F}_{n,i-1})\overset{p}{\to}\sigma^2$
\end{enumerate} where $\sigma^2$ is a constant. Then $S_{nk_{n}}=\sum_{i}Y_{ni}\overset{d}{\to}Z$, where the R.V. $Z$ has characteristic function $\exp(-\sigma^2t^2/2)$.
\end{thma}

Now, for $i>0$ we define $D_{i}=\hat{g}_{i}(X_i)-\Pi_{i-1}\hat{g}_{i}(X_{i-1})-c_{i}\{f(X_i)-\Pi_{i-1}f(X_{i-1})\}$, and take $k_n=n$, $\mathscr{F}_{ni}=\sigma(X_0,...,X_i)$, and $S_{ni}=n^{-1/2}\sum_{j=1}^{i}D_j$. From these definitions, we have $\mathscr{F}_{ni}\subset\mathscr{F}_{n,i+1}$ for $1\leq i<n$. 
We will verify Conditions~\ref{cond:condLind} and~\ref{cond:convVar} of Theorem~\ref{thm:mtgCLT} hold for $\{S_{ni},\mathscr{F}_{ni},1\leq i\leq k_n,n\geq 1\}$ defined in this way, following Section 17.4.2 of \citet{meyn2009markov}. In order to motivate this approach, we note that $S_{nn}=n^{-1/2}U_{n}$. 


Now, for $k=1,...,K$, we define $r_k(i)=k+(i-1)K$. For $t\geq k$, we define $m_k(t)=\max\{i\in\mathbb{N}:r_k(i)\leq t\}$. For checking the conditional Lindeberg condition~\ref{cond:condLind}, it is enough to show that \begin{align*}
	\sum_{i=1}^{m_k(n)}E\{Y_{n,r_k(i)}^2I(|Y_{n,r_k(i)}|>\epsilon)|\mathscr{F}_{n,r_k(i)-1}\}\overset{p}{\to}0
\end{align*} as $n\to\infty$ for each $k=1,...,K$. Conditions~\ref{cond:piStationary}--\ref{cond:drift} imply that for $k=1,...,K$, the subchains $(X_{k+Kt-1})_{t=1}^{\infty}$ are Harris recurrent with stationary measure $\pi$. Therefore, the Law of Large Numbers (Theorem 17.3.2 of \citet{meyn2009markov}) holds for each subchain. Consider an arbitrary $k$. For $i\geq 1,n\geq r_k(i)$, we have

\begin{align*}
	E\{D_{r_k(i)}^2I(|D_{r_k(i)}|>b)|\mathscr{F}_{n,r_k(i)-1}\}=h_k^b(X_{r_k(i)-1})
\end{align*} a.e. $P_\pi$ for some $\pi$-integrable function $h_k^b:X\to\mathbb{R}$. Therefore \begin{align*}
	&\underset{n}{\lim\sup}	\;\sum_{i=1}^{m_k(n)}E\{Y_{n,r_k(i)}^2I(|Y_{n,r_k(i)}|>b)|\mathscr{F}_{n,r_k(i)-1}\}\\
	&=\underset{n}{\lim\sup}\;n^{-1}\sum_{i=1}^{m_k(n)}E\{D_{r_k(i)}^2I(|D_{r_k(i)}|>n^{1/2}b)|\mathscr{F}_{n,r_k(i)-1}\}\\
	&\leq \underset{n}{\lim\sup}\;n^{-1}\sum_{i=1}^{m_k(n)}E\{D_{r_k(i)}^2I(|D_{r_k(i)}|>b^*)|\mathscr{F}_{n,r_k(i)-1}\}\\
	&\leq K^{-1}\underset{n\to\infty}{\lim\sup}	\{m_k(n)-1)\}^{-1}\sum_{t=1}^{m_k(n)}h_k^{b^*}(X_{r_k(i)-1})\\
	&=K^{-1}\int\pi(dx)h_k^{b^*}(x)
\end{align*} a.e. $P_\pi$ for any $b^*>0$, where the first equality follows from the definition of $Y_{ni}$, and the last equality follows from applying the Law of Large numbers to the subchain $\{X_{r_k(i)-1}\}_{i=1}^{\infty}$. Now, from the properties of conditional expectation, and the dominated convergence theorem, we can find a sequence $b_j\uparrow\infty$ for which $\int\pi(dx)h_k^{b_j}(x)\leq j^{-1}$ for each $j$. Thus, we obtain \begin{align*}
	&\underset{n}{\lim\sup}	\;\sum_{i=1}^{m_k(n)}E\{Y_{n,r_k(i)}^2I(|Y_{n,r_k(i)}|>b)|\mathscr{F}_{n,r_k(i)-1}\}\leq (jK)^{-1}
\end{align*} almost surely for each $j$, so the event \begin{align*}
	&\left\{\underset{n}{\lim\sup}	\;\sum_{i=1}^{m_k(n)}E\{Y_{n,r_k(i)}^2I(|Y_{n,r_k(i)}|>b)|\mathscr{F}_{n,r_k(i)-1}\}=0\right\}\\
	&=\cap_j\left\{\underset{n}{\lim\sup}	\;\sum_{i=1}^{m_k(n)}E\{Y_{n,r_k(i)}^2I(|Y_{n,r_k(i)}|>b)|\mathscr{F}_{n,r_k(i)-1}\}\leq (jK)^{-1}\right\}
\end{align*} has probability 1. Repeating this argument for each $k=1,...,K$ verifies the conditional Lindeberg condition~\ref{cond:condLind}. 

To verify the variance convergence in condition~\ref{cond:convVar}, we use the Law of Large Numbers on each subchain again to obtain
$\sum_{i}E(Y_{ni}^2|\mathscr{F}_{n,i-1})\overset{a.s.-P_{\pi}}{\to}\sigma^2$ where
\begin{align}
	&\sigma^2 =K^{-1}\sum_{k=1}^{K}\int \pi(dx)\Pi_k(x,dy)[\hat{g}_{\sigma(k)}(y)-\Pi_{k}\hat{g}_{\sigma(k)}(x)-c_{\sigma(k)}\{f(y)-\Pi_{k}f(x)\}]^2\nonumber\\
	&=K^{-1}\sum_{k=1}^{K}\left[\braket{\hat{g}_{\sigma(k)}-c_{\sigma(k)}f,\hat{g}_{\sigma(k)}-c_{\sigma(k)}f}\right.\nonumber\\&\left.-\braket{\Pi_k\hat{g}_{\sigma(k)}-c_{\sigma(k)}\Pi_kf,\Pi_k\hat{g}_{\sigma(k)}-c_{\sigma(k)}\Pi_kf}\right]\label{eq:sigma2}
\end{align} The convergence in probability in Condition \ref{cond:convVar} of Theorem~\ref{thm:mtgCLT} then follows immediately from the almost sure convergence. Thus by Theorem~\ref{thm:mtgCLT}, we have $S_{nn}\overset{d}{\to}Z$ where $Z$ has characteristic function $\exp(-\sigma^2t^2/2)$.

We now deal with the remainder term $\hat{g}_0(X_0)-\hat{g}_M(X_M)-c_0f(X_0)+c_Mf(X_M)$. Clearly, $M^{-1/2}\hat{g}_0(X_0)-c_0f(X_0)\overset{a.s.}{\to} 0$ as $M\to\infty$. Additionally, from the Law of Large Numbers applied to each subchain, \begin{align*}
	\sum_{t=0}^{M-1}M^{-1}\{\hat{g}_k(X_{k+Kt-1})-c_kf(X_{k+Kt-1})\}^2{\to}\int \pi(dx)\{\hat{g}_k(x)-c_kf(x)\}^2<\infty\end{align*} almost surely as $M\to\infty$ for each $k=1,...,K$. Therefore, $M^{-1/2}\{\hat{g}_M(X_M)+c_Mf(X_M)\}\overset{a.s.}{\to}0$ also.

Applying Slutsky's Theorem, we obtain $M^{-1/2}S_M\overset{d}{\to} Z$ where $Z$ has characteristic function $\exp(-\sigma^2t^2/2)$.

Now, we have $\sigma^2=K^{-1}\sum_{k=1}^{K}B_k$ where
\begin{align*}
	&B_k=\braket{\hat{g}_{\sigma(k)}-c_{\sigma(k)}f,\hat{g}_{\sigma(k)}-c_{\sigma(k)}f}-\braket{\Pi_k\hat{g}_{\sigma(k)}-c_{\sigma(k)}\Pi_kf,\Pi_k\hat{g}_{\sigma(k)}-c_{\sigma(k)}\Pi_kf}\\
	&=\braket{\hat{g}_{\sigma(k)}+\Pi_k\hat{g}_{\sigma(k)},\hat{g}_{\sigma(k)}-\Pi_k\hat{g}_{\sigma(k)}}-2c_{\sigma(k)}(\braket{f,\hat{g}_{\sigma(k)}}-\braket{\Pi_kf,\Pi_k\hat{g}_{\sigma(k)}})\\
	&+c_{\sigma(k)}^2(\braket{f,f}-\braket{\Pi_kf,\Pi_kf})\;\;\;\;\;(k=1,...,K).
\end{align*} Note \begin{align*}
	&\sum_{k=1}^{K}\braket{\hat{g}_{\sigma(k)}+\Pi_k\hat{g}_{\sigma(k)},\hat{g}_{\sigma(k)}-\Pi_k\hat{g}_{\sigma(k)}}\\
	&=\sum_{k=1}^{K}\braket{\hat{g}_{\sigma(k)},\hat{g}_{\sigma(k)}}-\braket{\Pi_k\hat{g}_{\sigma(k)},\Pi_k\hat{g}_{\sigma(k)}}\\
	&=\sum_{k=1}^{K}\braket{\hat{g}_{k},\hat{g}_{k}}-\braket{\Pi_k\hat{g}_{\sigma(k)},\Pi_k\hat{g}_{\sigma(k)}}\\
	&=\sum_{k=1}^{K}\braket{\hat{g}_{k}+\Pi_k\hat{g}_{k},\hat{g}_k-\Pi_k\hat{g}_{\sigma(k)}}\\
	&=\sum_{k=1}^{K}\braket{g,g}+2\sum_{t=1}^{\infty}\braket{g,P_k^tg}
\end{align*} where the last equality used Proposition~\ref{prop:poisEq} to simplify $\hat{g}_k-\Pi_k\hat{g}_{\sigma(k)}$. Thus, \begin{flalign*}
&\sigma^2=K^{-1}\sum_{k=1}^{K}B_k=\braket{g,g}+2K^{-1}\sum_{k=1}^{K}\sum_{t=1}^{\infty}\braket{g,P_k^tg}\\
&+K^{-1}\sum_{k=1}^{K}c_{\sigma(k)}^2(\braket{f,f}-\braket{\Pi_kf,\Pi_kf})-2c_{\sigma(k)}(\braket{f,\hat{g}_{\sigma(k)}}-\braket{\Pi_kf,\Pi_k\hat{g}_{\sigma(k)}})
\end{flalign*}



We now extend to the multivariate case by the Cramer-Wold device. Let $f:X\to\mathbb{R}^d$, $g:X\to\mathbb{R}^p$, $C_{k}\in\mathbb{R}^{p\times d}$, and \begin{align*}
	S_M=\sum_{t=0}^{M-1}g(X_t)-C_t^Tf(X_t)+C_{t+1}^T\Pi_tf(X_t).
\end{align*}
$a\in\mathbb{R}^d$. Define $U_k=\int\pi(dx)\{ff^T-(\Pi_kf)(\Pi_kf^T)\}$ and $V_k=\int \pi(dx)\{f\hat{g}_{\sigma(k)}^T-(\Pi_kf)(\Pi_k\hat{g}_{\sigma(k)})^T\}$. Then we have $a^TM^{-1/2}S_M\overset{d}{\to}Z$ where $Z$ is a random variable with characteristic function $\exp(-a^T\Sigma_C at^2/2)$, with \begin{align*}
	&\Sigma_C=\int\pi(dx)gg^T+2K^{-1}\sum_{k=1}^{K}\sum_{t=1}^{\infty}\int\pi(dx)g(P_kg)^T\\
	&+K^{-1}\sum_{k=1}^{K}C_{\sigma(k)}^TU_kC_{\sigma(k)}-C_{\sigma(k)}^TV_k-V_k^TC_{\sigma(k)}.
\end{align*} Since this holds for arbitrary $a$, we have by the Cramer-Wold Theorem that $M^{-1/2}S_M\overset{d}{\to}Z$, where $Z$ is a random variable with characteristic function $\exp(-t^T\Sigma_Ct/2)$.

Finally, we extend from the multivariate case with initial measure $\pi$, to the multivariate case with initial measure $\nu\neq \pi$. In this case, the desired convergence in distribution can be shown to hold via a coupling argument, as in \citet{roberts2004general}. We sketch the proof here. We construct on the same probability space two Markov chains $\{X_t\}_{t=0}^{\infty}$ and $\{\tilde{X}_t\}_{t=0}^{\infty}$, with initial law $\nu\times \pi$ for $(X_0,\tilde{X}_0)$. Then, we update the chains using a joint transition kernel chosen so that \begin{enumerate}
	\item each chain is marginally a Markov chain with transition kernel $\Pi_t=\Pi_{\sigma^t(1)}$ at time $t$, and
	\item $X_t=\tilde{X}_t$ for all $t>t_0$, for some random $t_0$, almost surely.
\end{enumerate} The aperiodicity assumption~\ref{cond:aperiodicity} and the geometric drift to the petite set $C$ in Assumption~\ref{cond:drift} ensure such a transition kernel can be constructed.

Then, $M^{-1/2}(S_M-\tilde{S}_M)\overset{a.s.}{\to}0$, where $S_M=\sum_{t=0}^{M-1}g(X_t)-C_t^T\{f(X_t)-\Pi_{t}f(X_t)\}$ and $\tilde{S}_M=\sum_{t=0}^{M-1}g(\tilde{X}_t)-C_t^T\{f(\tilde{X}_t)-\Pi_{t}f(\tilde{X}_t)\}$. Thus, from Slutsky's theorem, $M^{-1/2}S_M\overset{d}{\to}Z$ where $Z$ has characteristic function $\exp(-t^T\Sigma_Ct/2)$.

Now, we show that $\Sigma_C$ is minimized when $C_{\sigma(k)}=U_k^{\dagger}V_k$. First, we show $U_ka=0_{1\times d}$ implies $a^TV_k=0$. To see this, note that $U_ka=0$ implies \begin{align*}
	\int \pi(dx)\Pi_k(x,dy)a^T\{f(y)-\Pi_kf(x)\}\{f(y)-\Pi_kf(x)\}^Ta=0
\end{align*} so that $a^T\{f(y)-\Pi_kf(x)\}=0$ a.e. $\lambda_k$, where $\lambda_k$ is the measure on $(X^2,\mathscr{F}^2)$ defined by $\lambda_k(A\times B)=\int \pi(dx)\Pi_k(x,dy)I(x\in A,y\in B)$. In this case,
\begin{align*}
	&a^TV_k=\int\pi(dx)a^T\{f\hat{g}_{\sigma(k)}-\Pi_kf\Pi_k\hat{g}_{\sigma(k)}^T\}\\&=\int \lambda_k(dx\times dy)a^T\{f(y)-\Pi_kf(x)\}\{\hat{g}_{\sigma(k)}(y)-\Pi_k\hat{g}_{\sigma(k)}(x)\}\\
	&=0_{1\times d}.
\end{align*}

Finally, we note that $\Sigma_C$ depends on $C_{\sigma(k)}$ only through the term $K^{-1}(C_{\sigma(k)}^TU_kC_{\sigma(k)}-C^T_{\sigma(k)}V_k-V_k^TC_{\sigma(k)})$. By Lemma~\ref{lem:completeSquare}, this term is minimized when $C_{\sigma(k)}=U_k^{\dagger}V_k$. This completes the proof.\qed

\textbf{Proof of Corollary~\ref{cor:aVar0}}:
First, we obtain the simplified expression $V=K^{-1}\sum_{k=1}^{K}V_k=\int\pi(dx)fg^T$ for $V$. Under the Gibbs kernel assumption~\ref{cond:projection},

\begin{align*}
	&V=K^{-1}\sum_{k=1}^{K}V_k=K^{-1}\sum_{k=1}^{K}\int\pi(dx)\{f\hat{g}_{\sigma(k)}^T-\Pi_kf(\Pi_k\hat{g}_{\sigma(k)}^T)\}\\
	&=K^{-1}\sum_{k=1}^{K}\int\pi(dx)\{f\hat{g}_{k}^T-\Pi_kf(\Pi_k\hat{g}_{\sigma(k)}^T)\}\\
	&=K^{-1}\sum_{k=1}^{K}\int\pi(dx)\{f\hat{g}_{k}^T-f(\Pi_k\hat{g}_{\sigma(k)}^T)\}\\
	&=K^{-1}\sum_{k=1}^{K}\int\pi(dx)fg^T=\int \pi(dx)fg^T
\end{align*} where the second line rearranged the sum of the $f\hat{g}_{\sigma(k)}$ terms, and the third line used the equality $\int\pi(dx)\Pi_kf(\Pi_k\hat{g}_{\sigma(k)})^T=\int\pi(dx)f(\Pi_k\hat{g}_{\sigma(k)})^T$ from reversibility and idempotence of $\Pi_k$. The last line follows from Proposition~\ref{prop:poisEq}. 

In general, we have \begin{align*}
	&K^{-1}\sum_{k=1}^{K}C^TU_kC-V_k^TC-C^TV_K=C^TUC-C^TV-V^TC,
\end{align*}and \begin{align*}
	\Sigma_C&=\int \pi(dx)gg^T+K^{-1}\sum_{k=1}^{K}\sum_{t=1}^{\infty}\int\pi(dx)\{g(P_kg)^T+(P_kg)g^T\}\\
		&+K^{-1}\sum_{k=1}^{K}C_{\sigma(k)}^TU_kC_{\sigma(k)}-C_{\sigma(k)}^TV_k-V_k^TC_{\sigma(k)}.\\
		&=\int \pi(dx)gg^T+K^{-1}\sum_{k=1}^{K}\sum_{t=1}^{\infty}\int\pi(dx)\{g(P_kg)^T+(P_kg)g^T\}\\
		&+C^TUC-C^TV-V^TC
\end{align*} which is the representation of $\Sigma_C$ given in Corollary~\ref{cor:aVar0}.

Now, we show that $Ua=0$ implies $a^TV=0$, so that we may apply Lemma~\ref{lem:completeSquare} to the term $C^TUC-C^TV-V^TC$. First, we have that $Ua=0$ implies $a^TUa=0$, so from Lemma~\ref{lem:smallCondExp}, we have $Ua=0$ implies $a^Tf=b$ a.e. $\pi$. In this case $a^TV=0$ from the same reasoning as in the proof of Theorem~\ref{thm:vecG}. Thus, Lemma~\ref{lem:completeSquare} shows that $C^TUC-C^TV-V^TC$ is minimized when $C=\tilde{C}$, where $\tilde{C}=U^{\dagger}V$. Since $\Sigma_C$ depends on $C$ only through $C^TUC-C^TV-V^TC$, we have that $\Sigma_C$ is minimized at $C=\tilde{C}$. This completes the proof.\qed

\textbf{Proof of Corollary~\ref{cor:naiveVec}}:

We have \begin{align*}
	\Sigma_1&=\Sigma_0+K^{-1}\sum_{k=1}^{K}C_{\sigma(k)}^TU_kC_{\sigma(k)}-C_{\sigma(k)}^TV_k-V_k^TC_{\sigma(k)}\\
	&=\Sigma_0+K^{-1}\sum_{k=1}^{K}U_k-2V_k\\
	&=\Sigma_0-\int \pi(dx)gg^T-K^{-1}\sum_{k=1}^{K}\int \pi(dx)(\Pi_kg)(\Pi_kg)^T\leq \Sigma_0,
\end{align*} where the first equality used Theorem~\ref{thm:vecG}, and the second equality used $C_k=I_{d\times d}$ for each $k=1,...,K$ and the fact $f=g$. The third equality results from applying identity $K^{-1}\sum_{k=1}^{K}V_k=\int\pi(dx)=\int\pi(dx)gg^T$. The inequality holds since both integrals are of nonnegative functions, so that the subtracted integrands are nonnegative.\qed

\textbf{Proof of Theorem~\ref{thm:detVrev}:} First, consider the case where $g:X\to\mathbb{R}$ and $f:X\to\mathbb{R}$, and $C\in\mathbb{R}$, so that $S_M$ is a sum of scalar terms. 

From Lemma~\ref{lem:Qlem}, we have $\sum_{t=0}^{\infty}|Q^tg|$ is square integrable with respect to $\pi$. Thus, $\hat{h}=-Cf+\sum_{t=0}^{\infty}Q^tg$ is square integrable and satisfies the Poisson equation $\hat{h}-Q\hat{h}=h$ a.e. $\pi$. 

 It can be then be shown from the same martingale central limit theorem approach as in Theorem~\ref{thm:vecG} that for the random sweep chain, we have $M^{-1/2}S_M\to Z$ where $Z$ is a random variable with characteristic function $\exp(\Sigma_C^{rev}t^2/2)$, with $\Sigma_C^{rev}=\braket{\hat{h},\hat{h}}-\braket{Q\hat{h},Q\hat{h}}$. We now show that $\Sigma_C^{rev}$ is as given in the statement of Theorem~\ref{thm:detVrev}. We have
 \begin{align}
	\Sigma_C^{rev}&=\braket{\hat{h},\hat{h}}-\braket{Q\hat{h},Q\hat{h}}\nonumber\\
	&=\braket{-C(f+Qf)+g+2\sum_{t=1}^{\infty}Q^tg,g-C(f-Qf)}\nonumber\\
	&=\braket{h-2CQf+2\sum_{t=1}^{\infty}Q^tg,h}\label{eq:hVarStart}
\end{align} where the first equality follows from the martingale CLT, and the second inequality follows from the identify $a^2-b^2=(a+b)(a-b)$. Also, \begin{align}
	2\sum_{t=1}^{\infty}\braket{h,Q^th}&=2\sum_{t=1}^{\infty}\braket{h,Q^tg-CQ^t(f-Qf)}\nonumber\\
	&=-2C\sum_{t=1}^{\infty}\braket{h,Q^t(f-Qf)}+2\sum_{t=1}^{\infty}\braket{h,Q^tg}\nonumber\\
	&=-2C\braket{h,Qf}+2C\underset{t\to\infty}{\lim}\braket{h,Q^{t+1}f}+2\sum_{t=1}^{\infty}\braket{h,Q^tg}.\nonumber
\end{align} Now, we define $\|f\|=\braket{f,f}^{1/2}$. From Lemma 2 in~\citet{burkholder1961iterates}, we have since $Q$ is positive and self-adjoint that there exists an idempotent, self adjoint operator $\bar{Q}$ such that $\underset{t\to\infty}{\lim}\|\bar{Q}r-Q^tr\|=0$ for any function $r:X\to\mathbb{R}$ with $\braket{r,r}<\infty$. But for such a $\bar{Q}$, we have $Q\bar{Q}f=\bar{Q}f$ a.e. $\pi$, since \begin{align*}
	&\|Q\bar{Q}f-\bar{Q}f\|\leq \|Q\bar{Q}f-Q^tf\|+\|Q^tf-\bar{Q}f\|\\
	&\leq \|\bar{Q}f-Q^{t-1}f\|+\|Q^tf-\bar{Q}f\|
\end{align*} and $\underset{t\to\infty}{\lim}\|\bar{Q}-Q^{t-1}f\|+\|Q^tf-\bar{Q}f\|=0$. Since $Q\bar{Q}f=\bar{Q}f$ a.e. $\pi$, we have $\Pi_1\bar{Q}f=\bar{Q}f$ a.e. $\pi$ and $\Pi_2\bar{Q}f=\bar{Q}f$ a.e. $\pi$. Thus, from Lemma~\ref{lem:asConst}, $\bar{Q}f$ is constant a.e. $\pi$, so that $\underset{t\to\infty}{\lim}\braket{h,Q^{t+1}f}=\braket{h,\bar{Q}f}=0$, since $\int\pi(dx)h(x)=0$. Therefore, $2\sum_{t=1}^{\infty}\braket{h,Q^th}=-2C\braket{h,Qf}+2\sum_{t=1}^{\infty}\braket{h,Q^tg}$, so that we may rewrite~\eqref{eq:hVarStart} as \begin{align*}
	\Sigma_C^{rev}=\braket{h,h}+2\sum_{t=1}^{\infty}\braket{h,Q^th}
\end{align*} 
where we used square integrability of the sum of $|Q^tg|$ in order to apply Fubini's Theorem in~\eqref{eq:hVarStart}. This verifies the representation for $\Sigma_C^{rev}$ in Theorem~\ref{thm:detVrev}.

We now show that the two representations~\eqref{eq:SigmaC} and~\eqref{eq:sigmaC_h} coincide under the assumptions of Theorem~\ref{thm:detVrev}. First, we observe \begin{align*}
	\sum_{t=1}^{\infty}\braket{h,Q^tg}&=\sum_{t=1}^{\infty}\braket{g,(P_1^t+P_2^t)g}-C\sum_{t=1}^{\infty}\braket{g,Q^t(f-Qf)}\\
	&=-C\braket{g,Qf}+C\underset{t\to\infty}{\lim}\braket{g,Q^{t+1}f}+\sum_{t=1}^{\infty}\braket{g,(P_1^t+P_2^t)g}\\
	&=-C\braket{g,Qf}+\sum_{t=1}^{\infty}\braket{g,(P_1^t+P_2^t)g}
\end{align*} Thus, \begin{align*}
	\braket{h,h}+\sum_{t=1}^{\infty}\braket{h,Q^th}&=\braket{h,h}-C\braket{h,Qf}-C\braket{g,Qf}+\sum_{t=1}^{\infty}\braket{g,(P_1^t+P_2^t)g}\\
	&=\braket{g,g}+\sum_{t=1}^{\infty}\braket{g,(P_1^t+P_2^t)g}\\
	&-2C\braket{f,g}+\sum_{k=1}^{2}C^2(\braket{f,f}-\braket{\Pi_kf,\Pi_kf})
\end{align*} which coincides with~\eqref{eq:SigmaC} from Theorem~\ref{thm:vecG}.

The extension to multivariate $g:X\to\mathbb{R}^d$, $f:X\to\mathbb{R}^p$, and $C:X\to\mathbb{R}^{d\times p}$ follows via the Cramer-Wold device. We have $M^{-1/2}S_M\to Z$ where $Z$ has characteristic function $\exp(t^T\Sigma_C^{rev}t/2)$ with \begin{align*}
	\Sigma_C^{rev}=\int\pi(dx)hh^T+2\sum_{t=1}^{\infty}\int\pi(dx)h(Q^th)^T
\end{align*} for the random sweep chain and \begin{align*}
	\Sigma_C=\int\pi(dx)hh^T+\sum_{t=1}^{\infty}\int\pi(dx)h(Q^th)^T.
\end{align*} for the deterministic sweep chain. 
The expression in Theorem~\ref{thm:detVrev} for the difference $\Sigma_{\tilde{C}}-\Sigma_{\bar{C}}$ is obtained by arithmetic. We observe that $\sum_{t=1}^{\infty}\int\pi(dx)h(Qh)^T$ is positive semidefinite since $Q$ is a positive, self-adjoint operator and therefore has a positive, self-adjoint square root $\tilde{Q}$ with $\tilde{Q}\tilde{Q}=Q$, so that $\int\pi(dx)h(Q^th)^T=\int\pi(dx)\tilde{Q}^th(\tilde{Q}^th)^T\geq 0$.

\qed

\textbf{Proof of Proposition~\ref{prop:poisEq}}: We first prove the result for univariate $g:X\to\mathbb{R}$. Define $\tilde{g}_k(x)=\sum_{t=0}^{\infty}P^{Kt}_{k}g(x)$, for each $k$. By Assumption~\ref{cond:aperiodicity}, the kernels $P_k^K$ are aperiodic. Additionally, from Assumption~\ref{cond:irreducibility}, Markov chains resulting from $P_k^{K}$ are irreducible. From the geometric drift condition \ref{cond:drift}, Theorem 15.0.1 of \citet{meyn2009markov} implies \begin{align}
	\sum_{t=0}^{\infty}|P_k^{Kt}g(x)|\leq RV_k(x)\label{eq:poisBound}
\end{align} for some $R<\infty$ and all $x\in X$, for $k=1,...,K$.

Recall the definition \begin{align*}
&\hat{g}_k(x)=\sum_{t=0}^{\infty}P_{k}^tg(x)\;\;\;\;\;\;k=1,...,K
\end{align*} We now show that $\sum_{t=0}^{\infty}|P_k^tg(x)|$ is square integrable with respect to $\pi$ for each $k$. First, we note that it is sufficient to prove $\sum_{t=0}^{\infty}|P_k^{m}\{P_{\sigma^m(k)}^{Kt}g\}|$ is square integrable with respect to $\pi$ for each $m=0,...,K-1$, since in this case \begin{align*}\sum_{t=0}^{\infty}|P_{k}^tg(x)|=\sum_{m=0}^{K-1}\sum_{t=0}^{\infty}|P_k^{m}\{P_{\sigma^m(k)}^{Kt}g\}|\end{align*}.

Now, we have \begin{align*}
	&\int \pi(dx) \left[P_{k}^{m}\left\{\sum_{t=0}^{\infty}|P_{\sigma^m(k)}^{Kt}g|\right\}\right]^2\leq \int \pi(dx) P_{k}^{m}\left[\left\{\sum_{t=0}^{\infty}|P_{\sigma^m(k)}^{Kt}g|\right\}^2\right]\\&\leq\int \pi(dx) P_{k}^{m}\{R^2V_{\sigma^m(k)}^2\}=\int \pi(dx) R^2V_{\sigma^m(k)}^2<\infty
\end{align*} for each $m=0,...,K-1$. The first inequality follows from Jensen's inequality, the second inequality follows from~\eqref{eq:poisBound} and the equality follows because the $\Pi_k$ preserve the stationary probability distribution $\pi$.

We may then apply Fubini's theorem for $\pi$ a.e. $x$ to obtain \begin{align*}
		&\int \pi(dx)\left[P_{k}^{m}\left\{\sum_{t=0}^{\infty}|P_{\sigma^m(k)}^{Kt}g|\right\}\right]^2\\
		&=\int\pi(dx)\left[\sum_{t=0}^{\infty}P_k^{m}\{|P_{\sigma^m(k)}^{Kt}g|\}\right]^2<\infty.
	\end{align*}

	Now, from Jensen's inequality, we have \begin{align*}
		&\sum_{t=0}^{\infty}|P_k^{m}\{P_{\sigma^m(k)}^{Kt}g\}|\leq \sum_{t=0}^{\infty}P_k^{m}\{|P_{\sigma^m(k)}^{Kt}g|\}.
	\end{align*} so that $\sum_{t=0}^{\infty}|P_k^{m}\{P_{\sigma^m(k)}^{Kt}g\}|$ is square integrable with respect to $\pi$, for each $m=0,...,K-1$. 

	Thus, $\sum_{t=0}^{\infty}|P_k^tg|$ is square integrable with respect to $\pi$, and also $\hat{g}_k$ is square integrable with respect to $\pi$.


Now, we verify $\hat{g}_k-\Pi_k\hat{g}_{\sigma(k)}=g$ a.e. $\pi$ for each $k$. Since $\Pi_k\sum_{t=0}^{\infty}|P_{\sigma(k)}^tg|$ is square integrable with respect to $\pi$, we have from Fubini's theorem that \begin{align*}
	\Pi_k\hat{g}_{\sigma(k)}=\Pi_k\sum_{t=0}^{\infty}P_{\sigma(k)}^tg=\sum_{t=0}^{\infty}\Pi_kP_{\sigma(k)}^tg=\sum_{t=1}^{\infty}P_k^tg
\end{align*} for $\pi$ a.e. $x$, so that $\hat{g}_k-\Pi_k\hat{g}_{\sigma(k)}=g$ for $\pi$ a.e. $x$.

Now, for general $g:X\to\mathbb{R}^d$, we have $|a^Tg|\leq V_k$ from Assumption~\ref{cond:smallg} whenever $\|a\|_2\leq 1$. In particular, taking $a$ to be the vectors $e_i,i=1,...,d$ with $e_i$ having 1 in the $i$th position and $0$ elsewhere, we see from the previous reasoning that the conclusions of the Proposition still hold. We have $\hat{g}_k-\Pi_k\hat{g}_{\sigma(k)}=g$ a.e. $\pi$. Additionally $\int \pi(dx)\{\sum_{t=0}^{\infty}|P_k^tg|\}^T\{\sum_{t=0}^{\infty}|P_k^tg|\}<\infty$, so that the sum $\sum_{t=0}^{\infty}P_k^tg$ converges absolutely, elementwise, for $\pi$ a.e. $x$, and each of the $d$ components of $\sum_{t=0}^{\infty}|P_k^tg|$ are square integrable with respect to $\pi$.\qed

\textbf{Proof of Proposition~\ref{prop:lwkProp}}: First, we show that the LWK conditioning approach is, to within an asymptotically negligible term, an instance of the control variate scheme in~\eqref{eq:cv0} with $C=2I_{d\times d}$. Note that $Qg=(\Pi_1g+\Pi_2g)/2=g/2+\Pi_1g/2$. Define $H_M=\sum_{t=0}^{M-1}\Pi_1g(X_t)=\sum_{t=0}^{M-1}g(X_t)-2\{g(X_t)-Qg(X_t)\}$. Then for $S_M=\sum_{t=0}^{M-1}g(X_t)-2\{g(X_t)-\Pi_{\sigma^t(1)}g(X_t)\}$, we have from Lemma~\ref{lem:Qlem1} that $M^{-1/2}(S_M-H_M)\to 0$, so that $M^{-1/2}H_M$ has the same asymptotic distribution as $M^{-1/2}S_M$. Thus, $\Sigma_{\text{LWK}}=\Sigma_2$.

Since $\Pi_2g=g$ a.e. $\pi$, we have $U=U_1/2=(A-B)/2$. We also have $V=B$, and the term $C^TUC-V^TC-C^TV$ in the representation of $\Sigma_{C}$ in Corollary~\ref{cor:aVar0} can be written as $C^T(A-B)C/2-C^TA-A^TC$. Now, $\tilde{C}=U^{\dagger}V=2(A-B)^{-1}A$. Substituting $C=\tilde{C}$ and $C=2I_{d\times d}$ into~\eqref{eq:sigma2} and subtracting yields \begin{align*}
		&\Sigma_{\tilde{C}}-\Sigma_{2}=-2A(A-B)^{-1}A-(-2A-2B)\\
		&=-2A^T(A-B)^{-1}A+2(A-B)(A-B)^{-1}A+2B(A-B)^{-1}(A-B)\\
		&=-2B(A-B)^{-1}B\leq 0,
\end{align*} where the first equality used the symmetry of $A$ and $B$, and the final inequality used the fact that $A-B$ is positive semidefinite. When $B$ is positive definite, the inequality is strict.

Similarly, \begin{align*}
	\Sigma_{2}-\Sigma_{1}=-2(B+A)-\{(A-B)/2-2A\}=-(A+3B)/2<0
\end{align*} since $A$ is positive definite from Assumption~\ref{cond:variance} and $B$ is positive semidefinite. Finally, $\Sigma_{1}-\Sigma_0=-(B+3A)/2<0$. This completes the proof.\qed

\putbib
\end{bibunit}
\end{document}